\documentclass[12pt]{article}
\usepackage{amsmath, amssymb, amsthm}
\usepackage{graphicx}
\date{}
\input epsf.sty
\usepackage{lineno}
\begin{document}

\begin{center}
\baselineskip .2in {\large\bf }
\end{center}

%\begin{center}
%\baselineskip .2in {\large\bf Bifurcation and chaos in a multi-delayed
%prey-predator system with habitat complexity}
%\end{center}
%\begin{center}
%\baselineskip .2in {\large\bf Combined effect of negative and positive feedback delays on a prey-predator system with prey refuge}
%\end{center}

\begin{center}
\baselineskip .2in {\large\bf Complex dynamics generated by negative and positive feedback delays of a prey-predator system with prey refuge: Hopf bifurcation to Chaos}
\end{center}

\begin{center}
\baselineskip .2in {\bf  Debaldev Jana$^a$, R. Gopal$^b$ and M. Lakshmanan$^{c}$}

{\small\it $^a$Department of Mathematics,}
\\ SRM University, \\ Kattankulathur - 603 203, Tamil Nadu, India \\
{\small\it $^b$School of Electrical and Electronics Engineering,}
\\ SASTRA University, \\ Thanjavur -613 401, Tamil Nadu, India \\
{\small\it $^c$Centre for Nonlinear Dynamics,}
\\ School of Physics, Bharathidasan University,\\
 Tiruchirapalli - 620 024, India.\\

\end{center}

%\begin{center}
%\baselineskip .2in {\bf}
%
%{\small\it }
%\\  \\
%
%\end{center}

\begin{abstract} Various field and laboratory experiments show that
prey refuge plays a significant role in the stability of
prey-predator dynamics. On the other hand, theoretical studies show
that delayed system exhibits a much more realistic dynamics than its
non-delayed counterpart. In this paper, we study a multi-delayed
prey-predator model with prey refuge.
We consider modified Holling Type II response function that incorporates the effect of prey refuge
and then introduce two discrete delays in the model system. A negative feedback
delay is considered in the logistic prey growth rate to represent
density dependent feedback mechanism and a positive feedback delay is considered to represent the gestation time of the predator. Our study reveals that the system exhibits different dynamical behaviors, viz., stable coexistence, periodic coexistence or
chaos depending on the values of the delay parameters and degree of prey refuge. The interplay between two delays for a fixed value of prey refuge has also been determined. It is noticed that these delays work in a complementary fashion. In addition, using the normal form theory
and center manifold argument, we derive the explicit formulae for determining the direction of the bifurcation, the stability
and other properties of the bifurcating periodic solutions.

\end{abstract}

%
%\noindent $^{*}$ { Corresponding Author; E-mail:}
%
%\noindent {$\dag$Research is supported by UGC (Dr. D. S. Kothari Postdoctoral Fellowship),
%India; No.F.4-2/2006(BSR)/13-1004/2013(BSR).}
%
%\newpage

%\textbf{A prey-predator system with prey refuge and one negative feedback delay in prey growth term and a positive feedback delay as gestation period in predator growth term are considered.
%System is stable below the critical values of the delay (in case of any one delay is running and other is absent) or delays (in case of both delays are running) and above which it is unstable. System bifurcates (Hopf) at the critical value/ values of the delay parameters. These results are clearly depicted by analytical as well as numerical results.
%The system exhibits different dynamical behaviors, viz., stable coexistence, periodic coexistence (period 1, 2, 4, 5, 6, 8) or even chaos depending up on the different parametric conditions of the delay parameters and degree of prey refuge. These dynamical aspects are clearly featured by extensive numerical simulation including a two dimensional phase plane with varying two delays and other parameters are keeping constant.  Regular to chaotic behavior of the system is depicted by bifurcation diagram by varying one parameter and chaos is tested by positive Lyapunov exponent test.
%Using the normal form theory and center manifold argument, we derive the explicit formulae for determining the direction of the bifurcation, the stability and other properties of the bifurcating periodic solutions.}

\section{Introduction}

Interaction between food and its eater is the basic rule of nature. Interrelationship \cite{A01,A84} largely depends and varies upon
the structure of habitat and more prominently their habitat selection \cite{TWA01,J06}. Refuge means a place or state of safety \cite{JR16}. The predation risk often induces the use of safer
'refuge' habitats by prey population and refuge use under
predation risk is commonly observed in a wide range of
systems \cite{LD90,BET91,L98,BK04,C05,SB05,C09}. Refuge habitats can
change the dynamic behavior of prey-predator interaction by decreasing the predation risk \cite{A01,TWA01,L73}.  Ray and Stra$\check{\mbox{s}}$kraba \cite{RS01} notice that the prey species, i.e., detritivorous fish and their predator species, carnivorous fish coexist in Sundarban
Mangrove ecosystem. In the pristine part of the ecosystem where forest is dense, the production of detritus is high and the
detritivorous fish can easily take refuge in the densely inundated bushy part of the forest to avoid predation by carnivorous fish. In this
part the production and densities of both prey and predator fish are high and coexist with each other. But in the reclaimed part of
the forest where there is huge anthropogenic stress, production of detritus is less and due to lack of bushy part of the forest the area
for refuge of prey species is minimum. In this reclaimed area prey density (detritivorous fish) is reduced at an alarming level but predator fish
population is also slightly reduced because this fish population has switched its food habit to other fish and animals \cite{RET08}.
Predator's functional response is assumed to be one of the most important components of prey-predator interaction because
predation intensity may change the shape of the community structure and
ecosystem properties \cite{SET85}. It is well established that prey refuge reduces the predation rate and the hypothesis is that there exists an inverse relationship between predation
rate and degree of prey refuge \cite{SS82}. Therefore, it is important to incorporate the effect of prey refuge in the predator response function when we study the prey-predator model. A generalized Gause-type prey-predator model is given by
$$\begin{array}{l}
\dot{x}(t)=x(t)f(x(t))-y(t)g(x(t)),\\
\dot{y}(t)=\theta y(t)g(x(t))-dy(t),\end{array}\eqno{(1.1)}$$
 where $x$ and $y$ denote the prey and predator densities, respectively. $d$ is the death rate constant and $\theta$ is the conversion efficiency of the predator. $f(x)$ is the specific growth rate of prey in the absence of predator and $g(x)$ is the predator response function.
 It is shown that the most commonly used Holling Type II functional response, given by $g(x(t))=\frac{\alpha x }{1+\alpha hx}$, can be modified to $g(x(t))=\frac{\alpha (1-m)x }{1+\alpha (1-m)hx}$ \cite{K05,J13,J14,JET15} in the presence of prey refuge. Here $\alpha$ is the prey attack coefficient, $h$ is the handling time and $m ~(0<m<1)$ is the degree or strength of prey refuge. For example, $m = 0.30$ implies that prey-predator interaction will be reduced by $30\%$ due to prey refuge. If $m=0$, i.e. in the absence of prey refuge, it will be exactly the Holling Type II response function.

Assume that the prey grows logistically to environmental carrying capacity $k$ in the absence of  predator with intrinsic growth rate $r$. In this case, the generalized Gause-type prey-predator model (1.1) with modified Type II response function that incorporates the effect of prey refuge takes the following form:
$$\begin{array}{l}
\dot{x}(t)=rx(1-\frac{x}{k})-\frac{\alpha (1-m)x y}{1+\alpha (1-m)hx},\\
\dot{y}(t)= y\bigg[\frac{\theta \alpha (1-m)x}{1+\alpha (1-m)hx}-d\bigg].
\end{array}\eqno{(1.2)}$$

Models with delay are much more realistic, as in reality
time delays occur in almost every biological situation \cite{M89}
and they are assumed to be one of the reasons for regular fluctuations in
population density \cite{M81,XR01}. Assuming that reproduction of predator after
consuming prey is not instantaneous but mediated by some time
lag required for gestation, we considered a delay in the predator's numerical response of the above model and determined the critical value of the delay parameter below which the system was stable and above which it was unstable \cite{BJ11,J14,JET14,JET15}.
Recently, many researchers have studied the prey-predator interaction with two delays \cite{J14,JET15,YC06,SET08,NET06,XET11,LET12,JET16}.
In \cite{YC06}, the authors study a Lotka-Volterra prey-predator system with two delays.
Considering the sum of the two delays as a bifurcation parameter, they show that the system undergoes a Hopf bifurcation when the total delay exceeds some critical value. They also determine the direction and stability of the bifurcating periodic solutions. Song et al. \cite{SET08} study exactly the same model with additional analysis of global existence of periodic solutions. Nakaoka et al. \cite{NET06} and Xu et al. \cite{XET11} study similar type of Lotka-Volterra prey-predator model with two delays. It is to be mentioned that these studies \cite{YC06,SET08,NET06,XET11} consider the predator's functional response as Holling Type I and treat total delay as the bifurcation parameter where as stability and bifurcation analysis of a diffusive prey-predator system in Holling type III functional response with prey refuge is studied by \cite{YW14}. Diffusive prey-predator system with constant prey refuge and delay is studied by \cite{YZ16} and on the same type system with hyperbolic mortality is studied by \cite{YZ16a}. A three-species prey-predator system with two delays is studied in \cite{LET12}.
By choosing the sum of two delays as a bifurcation parameter, the authors show that a Hopf bifurcation
at the positive equilibrium of the system can occur as the total delay crosses some critical values and then determine the
direction and stability of the bifurcating periodic solutions. Depending only on a single phase diagram, Nakoka et al. \cite{NET06} and Liao et al. \cite{LET12} comment that large delay may cause chaos.
To know the effect of multi-delays on the qualitative behavior of prey-predator system with prey refuge, we modify the system (1.2) with two discrete delays. One discrete delay $\tau_1$ is considered in the specific growth rate of prey to incorporate the effect of density dependent feedback mechanism which takes $\tau_1$ units of time to respond to changes in the prey population \cite{F87}. The second delay $\tau_2$ is considered in the predator response function and it is regarded as gestation period or reaction time of the predator \cite{K93}. We thus obtain the following multi-delayed prey-predator model in the presence of prey refuge as
$$\begin{array}{l}
\dot{x}(t) =rx\bigg(1-\frac{x(t-\tau_1)}{k}\bigg)-\frac{\alpha (1-m)x y}{1+\alpha (1-m)hx},\\
\dot{y}(t) = y\bigg[\frac{\theta \alpha (1-m)x(t-\tau_2)}{1+\alpha (1-m)hx(t-\tau_2)}-d\bigg].
\end{array}\eqno{(1.3)}$$
All parameters are assumed to be positive. The model system (1.3) has to be investigated
with initial conditions
\[x(\phi)=x_0=\phi_1(\phi)>0,~ y(\phi)=y_0=\phi_1(\phi)>0 ~\mbox{for}~ \phi \in [-\mbox{max}\{\tau_1, \tau_2\},0].\eqno{(1.4)}\]
%\begin{align}
%& x(\theta)=\phi_1(\theta)\geq0, y(\theta)=\phi_2(\theta)\geq0, \nonumber \\
%& \theta\in[-\mbox{max}\{\tau_1,\tau_2\},0],~~\phi_i(0)>0~~~(i=1,2),
%\end{align}
%where $\phi:[-\mbox{max}\{\tau_1,\tau_2\},0]\rightarrow \Re^2$ with norm\\
%$$||\phi||=\sup \limits_{-\tau\leq\theta\leq0}\{|\phi_1(\theta)|,|\phi_2(\theta)|\},$$
%such that $\phi=(\phi_1,\phi_2)$.\\
The objective of our paper is to study the prey-predator dynamics in the presence of prey refuge and two biological delays.
We address the question how the gestation delay $\tau_2$ of the predator affects individually and jointly with the delay $\tau_1$ in the logistic prey growth rate on the local stability of a prey-predator system with prey refuge. Another important issue that will be discussed is -- how the system behaves if the biological delays are large enough. Existence of chaos, if any, due to large delay will be investigated rigorously. The interplay between the degree of prey refuge and biological delays will also be investigated in detail.\\

The organization of the paper is as follows: In Section 2, we study the local stability of the system (1.3); direction and stability of Hopf bifurcation is studied in Section 3. Rigorous numerical simulations of the model system are performed in Section 4.
Finally, a summary is presented in Section 5.\\

\section{Stability Analysis}

\subsection{Positive Invariance}

Feasibility or biological positivity studies aim to objectively and rationally uncover the strengths and weaknesses of an existing or proposed model in the given environment. Therefore, it is important to show the positivity for the model system (1.3) as the system represents prey-predator
populations. Biologically, positivity ensures that the population never becomes negative and
it always survives. For proving this, we have the following theorem.\\

{\bf Theorem 2.1.}  {\it All the solution of (1.3) with initial conditions (1.4) are positive.}\\

{\bf Proof.} The model (1.3) can be written in the following form:\\
$$\mbox{Consider}~W=col(x,y)\in \Re_+^2,  (\phi_1(\theta),\phi_2(\theta))\in \it C_+ =([-\mbox{max}\{\tau_1,\tau_2\},0],\Re_+^2),$$
$$\phi_1(0),\phi_2(0)>0,$$
\begin{eqnarray}F(W)&=&\begin{pmatrix}
                        _{\displaystyle{F_1(W)}}\\
                        _{\displaystyle{F_2(W)}}\nonumber
                       \end{pmatrix}\nonumber\\
\nonumber\\
                    &=&\begin{pmatrix}
                        _{\displaystyle{x\Bigg[r\bigg(1-\frac{x(t-\tau_1)}{K}\bigg)-\frac{\alpha (1-m) y}{ay+(1-m)x}\Bigg]}}\\
                        _{\displaystyle{y\Bigg[\frac{\theta\alpha(1-m)x(t-\tau_2)}{1+\alpha(1-m)hx(t-\tau_2)}-d\Bigg]}}\nonumber
                       \end{pmatrix},\nonumber\\
\nonumber\end{eqnarray}\nonumber\\
model system (1.3) becomes
\[\dot W=F(W),\eqno{(1.3a)}\]
with $W(\phi)=  (\phi_1(\phi),\phi_2(\phi))\in \it C_+$ and $\phi_1(0),\phi_2(0)>0$. It is easy to check in system
 (1.3a) that whenever choosing $W(\phi) \in \Re_+$ such that $x=y=0,$ then\\
 $$F_i(W)\mid_{w_i=0,W \in \Re_+^2} \geq 0,$$
 with $w_1(t)=x(t), w_2(t)=y(t)$. Using the lemma given in \cite{XY96}, any solution of (1.3a) with $W(\theta)\in C_+,$ say $W(t)=W(t,W(\theta))$,
 is such that $W(t)\in\Re_+^2$ for all $t\geq 0.$ Hence the solution of the system (1.3a) exists in the region $\Re_+^2$ and all solutions
 remain non-negative for all $t>0$. Therefore, the positive orthant $\Re_+^2$ is an invariant region.

 \subsection{Uniform persistence}
In ecology, one question arises that the determining conditions which assure that the population abundance with time i.e., the solutions of the corresponding ecological system which are initially strictly positive do not approach the boundary of the cone as time evolves. Generally speaking, the uniformly persistent systems are those in which strictly positive solutions do not approach the boundary of the non-negative phase-space ($\Re^2_+$) as $t \rightarrow \infty$ \cite{BET86}. In other words, permanently coexistence (uniform persistence) implies the existence of a region in the phase space at a non-zero distance from the boundary, in which all the population vectors must lie ultimately. The uniform persistence is also most suitable from the point of applications since it rules out the possibility of one of the populations becoming arbitrarily close to zero and hence the risk of extinction due to small perturbation due to stochastic effects. The  concept of uniform persistence has been discussed by many researchers \cite{FW84,FS85,BC01}. In the next theorem, we prove the uniform persistence of the the model system (1.3) by means of the well known average Lyapunov function \cite{GH79}.\\

{\bf Theorem 2.2.}{\it The model system (1.3) is uniformly persistent if all the non-interior equilibrium points exist together with $(i)~ \frac{\beta_1}{\beta_2}>\frac{d}{r},~ (ii)~ m<1-\frac{d}{\alpha k(\theta-hd)},~ \theta>hd+\frac{d}{\alpha k}$.}

{\bf Proof} The system (1.3) has two boundary equilibrium points, $(i) E_0(0,0)$ and $(ii) E_1(k,0)$. For $(x,y) \in \Re^2_+,$ consider the following average Lyapunov function.
\[\rho(x,y,z) = x^{\beta_1}y^{\beta_2},\]
where $\beta_i>0,$ $i = 1,2$ are constants. It is easy to observe that the function $\rho(.)$ is a non-negative continuous function. Calculating the logarithmic derivative of $\rho(.)$ along solution of the model system (1.3), we obtain
\begin{eqnarray}
\nonumber
\Upsilon(x,y,z) &=& \frac{\dot{\rho}}{\rho} = \frac{\beta_1\dot{x}}{x}+\frac{\beta_2\dot{y}}{y}\\
\nonumber
&=& \beta_1\Big[r\bigg(1-\frac{x(t-\tau_1)}{k}-\frac{\alpha(1-m)y}{1+\alpha(1-m)hx}\bigg) \Big]+\beta_2\Big[\frac{\theta\alpha(1-m)x(t-\tau_2)}{1+\alpha(1-m)hx(t-\tau_2)}-d\Big].
\nonumber
\end{eqnarray}
To prove the system (1.3) to be uniformly persistent, we need to show that the above logarithmic derivative is positive at all the boundary equilibrium points $E_0$ and $E_1$ for some suitable choices of positive $\beta_i$'s. This condition is ensured at the origin and axial equilibria $E_0$ and $E_1$ by the choice $\frac{\beta_1}{\beta_2}>\frac{d}{r},$ $m<1-\frac{d}{\alpha k(\theta-hd)}$ and $\theta>hd+\frac{d}{\alpha k}$ respectively.

 \subsection{Local stability}

Ecological stability can refer to types of stability in a continuum ranging from regeneration via resilience (returning quickly to a previous state), to constancy to persistence. The precise definition depends on the ecosystem in question, the variable or variables of interest, and the overall context. In the context of conservation ecology, stable populations are often defined as ones that do not go extinct. Researchers applying mathematical models from system dynamics usually use Lyapunov stability. Local stability indicates that a system is stable over small short-lived disturbances. In ecology,
stress is given on the stability of coexistence of equilibrium points. We,
therefore, concentrate on the study of the interior equilibrium point of the system (1.3).
This system has only one interior equilibrium point given by
$E^{*}(x^{*},y^{*})$, where
$x^{*}= \frac{d}{\alpha (1-m) (\theta-h d)}$ and $y^{*}=\frac{r (k-x^*)
\{1+\alpha h (1-m) x^*\}}{\alpha k (1-m) }$. The equilibrium point will
be feasible if
$$\begin{array}{l}
    (i)~ m<1- \frac{d}{\alpha k (\theta-h d)} ~\mbox{ and}\\
    (ii)~\theta >hd +\frac{d}{\alpha k}.
  \end{array}$$
%[(i) and (ii) are the existence condition of $E^{*}(x^{*},y^{*})$, putting $x^*>0, y^*>0$, we will get (i) and (ii).]
Linearizing the system $(1.3)$ at $(x^*,y^*)$, we get
%[By using the linear transformation $X=x-x^*,Y=y-y^*$, putting these transformations in equation (1.3) and then for simplicity change new state variables $X$ and $Y$ by $x$ and $y$, we get (2.1)]
$$\begin{array}{l}
  \dot{x}(t)=\frac{\alpha^2(1-m)^2hx^{*}y^{*}}{\{1+\alpha (1-m)h x^{*}\}^2}x(t)-\frac{\alpha (1-m)x^*}{1+\alpha (1-m)h x^{*}}y(t)-\frac{rx^*}{k}x(t-\tau_1), \\
    \dot{y}(t)=\frac{\theta \alpha (1-m)y^*}{\{1+\alpha (1-m)h x^{*}\}^2}x(t-\tau_2).
  \end{array}\eqno{(2.1)}$$
The characteristic equation of the corresponding variational matrix is given by
$$\begin{array}{l}
    \lambda^2+A \lambda+B\lambda e^{-\lambda \tau_1}+C e^{-\lambda \tau_2}=0,
  \end{array}\eqno{(2.2)}$$
where
$$\begin{array}{lll}
      A=-\frac{\alpha ^2(1-m)^2 h x^* y^*}{\{1+\alpha (1-m)h x^*\}^2} ~(<0), &
      B=\frac{rx^*}{k} ~(>0)~~\mbox{and} & C=\frac{\theta \alpha^2(1-m)^2x^*y^*}{\{1+\alpha (1-m)h x^*\}^3} ~(>0). \\
    \end{array}  $$
One can consider several subcases. They are as follows,
\subsection{Case I: $\tau_1=0 ~\mbox{and}~ \tau_2=0$}
In the absence of all delays,
the characteristic equation $(2.2)$ becomes
$$\begin{array}{l}
    \lambda^2+(A+B)\lambda+C=0.
  \end{array}\eqno{(2.3)}$$
All roots of the equation $(2.3)$ will have negative real parts and the corresponding non-delayed system of the delay-induced system (1.3) will be locally asymptotically
stable around $E^{*}$ if and only if
$$(H_{1})~~~A+B>0 ~~\mbox{and} ~~C>0.$$
Note that $C$ is always positive whenever $E^{*}$ exists and $A+B>0$ if
$$\begin{array}{l}
    m>1- \frac{\theta+hd}{\alpha
kh(\theta-hd)}\mbox{ with} ~\theta >\frac{hd (1+\alpha
kh)}{\alpha kh-1}\mbox{ and}~ \alpha >\frac{1}{kh}.
  \end{array}$$
Thus, we state the following lemma for the stability of the non-delayed system.\\

\noindent {\bf Lemma 2.1.} The system (1.3) is locally asymptotically stable around $E^{*}$ in the absence of delay if
$$\begin{array}{l}
    (i)~\alpha >\frac{1}{kh},\\
    (ii)~\theta>\mbox{max}\big[\frac{hd(\alpha kh+1)}{\alpha
kh-1},hd+\frac{d}{\alpha k}\big] ~ \mbox{and}\\
    (iii) ~1- \frac{\theta+hd}{\alpha
kh(\theta-hd)}<m<1- \frac{d}{\alpha k (\theta-h d)}.
  \end{array}$$

\subsection{Case II: $\tau_1= 0 ~\mbox{and}~ \tau_2 \neq 0$}
If $\tau_{1}=0$ and $\tau_{2}>0$ then the characteristic
equation $(2.2)$ becomes
$$\begin{array}{l}
\lambda^2+(A+B)\lambda+Ce^{-\lambda \tau_2}=0\end{array}\eqno{(2.4)}$$
Let $i \omega ~(\omega>0)$
be a root of the equation $(2.4)$. Then it follows that
$$\begin{array}{l}
   \omega^2=C \mbox{cos}\omega \tau_2, \\
    (A+B)\omega=C \mbox{sin} \omega \tau_2.
  \end{array}\eqno{(2.5)}$$
This leads to
$$\begin{array}{l}
    \omega^4+(A+B)^2\omega^2-C^2=0.
  \end{array}\eqno{(2.6)}$$
If Lemma 2.1 holds, then the equation $(2.6)$ has a unique positive root $\omega_0^2$. Substituting
$\omega_0^2$ into $(2.5)$, we have
$$\begin{array}{l}
    \tau_{2_{n}}=\frac{1}{\omega_0}\mbox{cos}^{-1}\bigg(\frac{\omega_0^2}{C}\bigg)+\frac{2n \pi}{\omega_0},~
    n=0,1,2,....
  \end{array}$$
%Let, $$(H_2)~~~A+B<0.$$
%If $(H_1)$ and $(H_2)$ hold, equation $(2.6)$ will have two positive roots $\omega_{+}^2$ and $\omega_{-}^2$.
%Substituting $\omega_{\pm}$ into equation $(2.5)$, we get
%$$\begin{array}{l}
%    \tau_{2_{m}}^{\pm}= \frac{1}{\omega_{\pm}}\mbox{cos}^{-1}\bigg(\frac{\omega_{\pm}^2}{C}\bigg)+\frac{2m \pi}{\omega_{\pm}},~
 %   m=0,1,2,....
%  \end{array}$$
%If $\lambda(\tau_2)$ be the root of equation $(2.4)$, satisfying $\mbox{Re} \lambda(\tau_{2_{n}})=0$
%(respectively $\mbox{Re} \lambda(\tau_{2_{m}}^{\pm})=0$) and $\mbox{Im}\lambda(\tau_{2_{n}})=\omega_0$
%(respectively $\mbox{Im}\lambda(\tau_{2_{m}}^{\pm})=\omega_{\pm})$, we get
%$$\begin{array}{l}
%    \bigg[\frac{d}{d\tau_2}\mbox{Re}(\lambda)\bigg]_{\tau_{2}=\tau_{2_0},\omega=\omega_0}=
%    \frac{\omega^4+C^2}{\omega^2C^2}>0.
%  \end{array}$$
%Similarly, we can show that
%$$\begin{array}{l}
%    \bigg[\frac{d}{d\tau_2}\mbox{Re}(\lambda)\bigg]_{\tau_{2}=\tau_{2_m}^{+},\omega=\omega_{+}}>0 ~~~\mbox{and}~~~
%    \bigg[\frac{d}{d\tau_2}\mbox{Re}(\lambda)\bigg]_{\tau_{2}=\tau_{2_m}^{-},\omega=\omega_{-}}<0.
%  \end{array}$$
If $\lambda(\tau_{2})$
be the root of $(2.4)$ satisfying $Re \lambda(\tau_{2_{n}})=0$
and $Im
\lambda(\tau_{2_{n}})=\omega_{0}$, we get
$$\begin{array}{l}
    \bigg[\frac{d}{d\tau_2}Re(\lambda)\bigg]_{\tau_2=\tau_{2_0},\omega=\omega_0}=\frac{2{\omega}^2+(A+B)^2}{C^2}>0.
  \end{array}
$$
From Corollary (2.4) in Ruan and Wei \cite{RW03}, we have the
following conclusions.\\

%\noindent {\bf Lemma 2.2.} Assume $\tau_{1}=0$ and $(H_{1})$ is satisfied.
%Then the following conclusion holds:

%(i) The equilibrium $E^{*}$
%is asymptotically stable for $\tau_{2}<\tau_{2_{0}}$ and unstable
%for $\tau_{2}>\tau_{2_{0}}$. Furthermore, the system (1.3)
%undergoes a Hopf bifurcation at $E^{*}$ when
%$\tau_{2}=\tau_{2_{0}}$.
%
%(ii) If $H_{2}$ holds, then there is a positive integer $m$ such
%that the equilibrium is stable when
%$\tau_{2}\in[0,\tau_{2_{0}}^{+})\cup(\tau_{2_{0}}^{-},\tau_{2_{1}}^{+})\cup.....
%\cup(\tau_{2_{p-1}}^{-},\tau_{2_{p}}^{+})$,
%and unstable when
%$\tau_{2}\in[\tau_{2_{0}}^{+},\tau_{2_{0}}^{-})\cup(\tau_{2_{1}}^{+},\tau_{2_{1}}^{-})\cup.....
%\cup(\tau_{2_{p-1}}^{+},\tau_{2_{p-1}}^{-})\cup(\tau_{2_{p}}^{+},\infty)$.
%Further more, the system $(1.3)$ undergoes a Hopf bifurcation at
%$E^{*}$ when $\tau_{2}=\tau_{2_{m}}^{+}$,
%m=0,1,2,....

\noindent {\bf Lemma 2.2.} Assume $\tau_{1}=0$ and conditions of Lemma 2.1 hold. Then the interior equilibrium $E^{*}$ of the system (1.3)
is asymptotically stable for $\tau_{2}<\tau_{2_{0}}$ and unstable
for $\tau_{2}>\tau_{2_{0}}$. Furthermore, the system (1.3)
undergoes a Hopf bifurcation at $E^{*}$ when
$\tau_{2}=\tau_{2_{0}}$.

\subsection{Case III: $\tau_{1},\tau_{2}\neq0$ $\&$ $\tau_{2}$ \mbox{is within its stable range}}
In this case, we consider equation $(2.2)$
with $\tau_{2}$ in its stable interval and regard $\tau_{1}$ as a
parameter. Without loss of generality, we consider that the system
$(2.2)$ satisfies conditions of the Lemma 2.1. Let $i \omega
(\omega>0)$ be a root of equation $(2.2)$ and thus we obtain
$$\begin{array}{l}
    \omega^4+\widetilde{A}\omega^2+C^2+2\widetilde{B}\mbox{cos}{\omega \tau_2}+2\widetilde{C}\mbox{sin}
    {\omega \tau_2}=0.
  \end{array}\eqno{(2.7)}$$
where,
$$\begin{array}{lll}
    \widetilde{A}=A^2-B^2, & \widetilde{B}=-C\omega^2 ~\mbox{and} & \widetilde{C}=-AC\omega.
  \end{array}$$
We define,
$$\begin{array}{l}
F(\omega)=\omega^4+\widetilde{A}\omega^2+C^2+2\widetilde{B}\mbox{cos}{\omega \tau_2}+2\widetilde{C}\mbox{sin}
    {\omega \tau_2}.\end{array}$$
Then it is easy to check that $F(0)=C^2>0$ and $F(\infty)=\infty$. We
can obtain that equation $(2.7)$ has finite positive roots $\omega_{1},
\omega_{2},....., \omega_{k}$. For every fixed  $\omega_{i}, i=1,
2, ...., k $, there exists a sequence $\{\tau_{1_{i}}^{j}\mid j=1,
2, .... \}$, where

$$\begin{array}{l}
\tau_{1_{i}}^{j}=\frac{1}{\omega_1}\mbox{cos}^{-1}\bigg[\frac{C~\mbox{sin}\omega_{1}\tau_{2}
-A\omega_{1}}{B\omega_{1}}\bigg]+\frac{2i\pi}{\omega_{1}},~~i=1,2,...,k,~j=1,
2, ....
 \end{array}$$
such that (2.7) holds. Let
$\tau_{1_{0}}=\min\{\tau_{1_{i}}^{j}\mid i=1, 2, ...., k;j=1, 2,
....\}$. When $\tau_{1}=\tau_{1_{0}}$, equation $(2.2)$ has a
pair of purely imaginary roots $\pm i \omega_{1}$ for $\tau_{2}
\in [0,\tau_{2_{0}})$.

In the following, we assume that
$$(H_2)~~~~~~~~~~~\left[\frac{d}{d \tau_{1}}(Re
\lambda)\right]_{\lambda=i \omega_{1}}\neq0.$$ Therefore,
by the general Hopf bifurcation theorem for FDEs, we have
the following result on the stability and bifurcation of the
system equation $(1.3)$.\\

\noindent {\bf Lemma 2.3.} For the model system $(1.3)$, suppose conditions of the Lemma 2.1 are satisfied when $\tau_{2}\in[0,\tau_{2_{0}})$. Then the equilibrium $E^{*}$ is locally asymptotically stable when
$\tau_{1}\in(0,\tau_{1_{0}})$ and unstable if $\tau_{1} > \tau_{1_{0}}.$ The system $(1.3)$ undergoes a
Hopf bifurcation at $E^{*}$ when $\tau_{1}=\tau_{1_{0}}$.
%$$\begin{array}{l}
%\tau_{1_{m}}=\frac{1}{\omega_1}\mbox{cos}^{-1}\bigg[\frac{C~\mbox{sin}\omega_{1}\tau_{2}-A\omega_{1}}{B\omega_{1}}\bigg]+\frac{2m\pi}{\omega_{1}},~~m=0,1,2,3,...
% \end{array}$$

\subsection{Case IV: $\tau_2= 0 ~\mbox{and}~ \tau_1 \neq 0$}

\noindent  For this choice of the delay parameters we summarize our results in the following theorem. The proof follows similar arguments as the Lemma 2.2. in Cases II.\\

\noindent{\bf Lemma 2.4.} Assume that $\tau_{2}=0, \tau_{1}\neq0$ and the condition of Lemma 2.1. hold.
Then the equilibrium $E^*$ is locally asymptotically stable for $\tau_{1}<\bar{\tau}_{1_{0}}$ and unstable
for $\tau_{1}>\bar{\tau}_{1_{0}}$. Furthermore, the system $(1.3)$ undergoes Hopf-bifurcation when $\tau_{1}=\bar{\tau}_{1_{0}}$, where

$$\begin{array}{l}
    \bar{\tau}_{1_{0}}=\frac{1}{\bar{\omega}_0}\mbox{cos}^{-1}(-\frac{A}{B}),
  \end{array}
 $$
 where, $\bar{\omega}_0$ is the unique positive root of ${\omega}^4+(A^2-B^2-2C){\omega}^2+C^2=0$ and
 \begin{eqnarray}\left[\frac{d ({\mbox Re} \lambda (\tau))}{d
\tau}\right]^{-1}_{\tau={\bar{\tau}}}
&=&\frac{M+2 {\bar{\omega}}^2}{B^2
{\bar{\omega}}^2}>0,\nonumber\end{eqnarray}
since $M=A^2-B^2-2C=[\frac{rd\{\alpha h(1-m)(k-2x^*)-1\}}{\alpha \theta k(1-m)}]^2>0$.

\subsection{Case V: $\tau_{1},\tau_{2}\neq0$ $\&$ $\tau_{1}$ \mbox{is within its stable range}}

For this choice of the delay parameters we summarize our results in the following theorem. The proof follows similar arguments as the Lemma 2.3. in Cases III.\\

\noindent {\bf Lemma 2.5.} For the model system $(1.3)$, suppose conditions of the Lemma 2.1 are satisfied and $\tau_{1}\in[0,\bar{\tau}_{1_{0}})$. Then the equilibrium $E^{*}$ is locally asymptotically stable when
$\tau_{2}\in(0,\bar{\tau}_{2_{0}})$ and unstable if $\tau_{2} > \bar{\tau}_{2_{0}}.$ The system $(1.3)$ undergoes a
Hopf bifurcation at $E^{*}$ when $\tau_{2}=\bar{\tau}_{2_{0}}$, where

$$\begin{array}{l}
    \bar{\tau}_{2_0}=\frac{1}{\bar{\omega}_1}\mbox{cos}^{-1}\bigg(\frac{\bar{\omega}_1^2-B\bar{\omega}_1 \mbox{sin}\bar{\omega}_1 \tau_1}{C}\bigg),
  \end{array}
$$
where, $\bar{\omega}_1>0$ is the positive root of
$$\begin{array}{l}
    \omega^4+\bar{A}\omega^2-C^2+a\bar{C}\mbox{sin}\omega\tau_1+\bar{D}\mbox{cos}\omega\tau_1=0,
  \end{array}
$$
with $\bar{A}=A^2+B^2,~\bar{B}=-B\omega^3$ and $\bar{C}=AB\omega^2.$

In the following, we also assume that
$$(H_3)~~~~~~~~~~~\left[\frac{d}{d \tau_{2}}(Re
\lambda)\right]_{\lambda=i \bar{\omega}_{1}}\neq0.$$

\section{Direction and stability of Hopf bifurcation point}
 In the previous section, we have obtained the sufficient conditions to guarantee
that the system (1.3) undergoes Hopf bifurcation at $E^{*}$ when $\tau_{1} =
\tau_{{1}_{0}}=\widetilde{\tau}_1$ and $\omega_1=\widetilde{\omega}$ (say) and $\tau_2$ is within its stability range. In this section, we study its bifurcation
properties. The method we use is based on the normal form and the center manifold theory presented in Hassard et al. \cite{HET81}.\\

\noindent Without loss of generality, we assume that $\widetilde{\tau_{2}} < \tau_{{2}_{0}}$,
where $\widetilde{\tau_{2}} \in (0, \tau_{{2}_{0}})$. Let $x_{1}=x-{x}^*,~
x_{2}=y-y^*$ and $\tau_1=\widetilde{\tau}_1+\mu$ where $\mu\in R$.
%We normalize the delay with the scaling $t\mapsto(t/\tau_{1})$.
Then equation (1.3) is transformed into an FDE in $C=C([-1,0],R^2)$ as
 $$\dot{x}(t)=L_{\mu}(x_{t})+f(\mu,x_{t}).\eqno{(4.1)}$$
Here $x(t)={(x_{1},x_{2})^{T}}\in R^{2}$, and
$L_{\mu}:C\rightarrow R,~f:R \times C \rightarrow R$ are given by
$$L_{\mu}(\phi)=(\widetilde{\tau}_1+\mu)
  B_{1}
\begin{pmatrix}
  _{\phi_{1} (0)} \\
  _{\phi_{2} (0)}
\end{pmatrix} + (\widetilde{\tau}_1+\mu)B_{2} \begin{pmatrix}
  _{\phi_{1} (-\widetilde{\tau}_1)} \\
  _{\phi_{2} (-\widetilde{\tau}_1)}
\end{pmatrix}+
 (\widetilde{\tau}_1+\mu) B_{3}
\begin{pmatrix}
  _{\phi_{1} (-\widetilde{\tau}_2)} \\
  _{\phi_{2} (-\widetilde{\tau}_2)}
\end{pmatrix} \eqno{(4.2)}$$
\noindent with $B_{1}=\begin{pmatrix}
  _{\frac{\alpha^2 (1-m)^2 h x^* y^*}{\{1+\alpha (1-m)h x^*\}^2}} & _{-\frac{\alpha (1-m)x^*}{1+\alpha (1-m)h x^*}} \\
  _{0} & _{0}
\end{pmatrix}$,
\noindent  $B_{2}=\begin{pmatrix}
  _{-\frac{r x^*}{k}} & _{0} \\
  _{0} & _{0}\end{pmatrix}$,
\noindent  $B_{3}=\begin{pmatrix}
    _{0} & _{0} \\
    _{\frac{\theta \alpha (1-m)y^*}{\{1+\alpha (1-m)h x^*\}^2}} & _{0} \
  \end{pmatrix}$
and
$$f(\mu,\phi)=(\widetilde{\tau}_1+\mu) \begin{pmatrix}
  _{-\frac{r}{k}\phi_1(0)\phi_1(-\widetilde{\tau}_1)-\frac{\alpha (1-m)\phi_1(0)\phi_2(0)}{1+\alpha (1-m) h \phi_1 (0)}} \\
  _{\frac{\theta \alpha (1-m)\phi_1(-\widetilde{\tau}_2)\phi_2(0)}{1+\alpha (1-m) h \phi_1 (-\widetilde{\tau}_2)}}
\end{pmatrix}.\eqno{(4.3)}$$
\noindent By Riesz representation theorem, there exits a function
$\eta(\vartheta,\mu)$ of bounded variation for $\vartheta \in [-1,0]$
such that
$$L_{\mu} \phi=\int_{-1}^0 d \eta(\vartheta,\mu) \phi(\vartheta) ~~~~ {\mbox for}~
\phi \in C. \eqno{(4.4)}$$\noindent In fact, we can choose
$$\eta(\vartheta,\mu)=\left\{\begin{array}{ll}\displaystyle (\widetilde{\tau}_1+\mu)B_1,&\vartheta = 0, \\
  (\widetilde{\tau}_1+\mu)B_2 \delta (\vartheta+\widetilde{\tau}_1),&\vartheta \in [-\widetilde{\tau}_1,0), \\
  -(\widetilde{\tau}_1+\mu)B_3 \delta(\vartheta+\widetilde{\tau}_2),&\vartheta \in [-\widetilde{\tau}_2,-\widetilde{\tau}_1). \\
  \end{array}\right.\eqno{(4.5)}$$
\noindent For $\phi \in C^{1} ([-1,0], R^{2})$, we define
$$A(\mu) \phi = \left\{\begin{array}{ll}\displaystyle{\frac{d \phi (\vartheta)}{d \vartheta}} ,& \vartheta \in  [-\widetilde{\tau}_2,0),\\
\int_{-1}^{0} d \eta(\mu,s) \phi (s), & \vartheta  = 0,
\end{array}\right.$$
 and
$$R(\mu) \phi = \left\{\begin{array}{ll}
  0, & \vartheta \in [-\widetilde{\tau}_2,0) \\
  f (\mu,\phi),& \vartheta=0.
\end{array}\right.$$
\noindent Then system (4.1) is equivalent to
$$\dot{x} (t) = A(\mu) x_{t} + R(\mu) x_{t} ,\eqno{(4.6)}$$\noindent
where $x_{t}(\vartheta) = x (t+\vartheta) {\mbox ~~~for} ~~~\vartheta \in
[-1,0].$

\noindent For $\psi \in C^{1}([0,-1], (R^{2})^*),$ define
$$A^{*} \psi (s) = \left\{\begin{array}{ll}
  -\displaystyle{\frac{d \psi (s)}{d s}}, & s \in (0,1], \\
  \int_{-1}^{0} d {\eta}^T (t,0) \psi (-t), & s=0,
\end{array}\right.$$
\noindent and a bilinear inner product
$$\langle \psi (s),\phi (\vartheta) \rangle = \overline{\psi} (0) \phi (0)-
\int_{-1}^{0} \int_{\xi = 0}^{\vartheta} \overline{\psi} (\xi -
\vartheta) d \eta (\vartheta) \phi (\xi) d \xi ,\eqno{(4.7)}$$
where $\eta (\vartheta)$ = $\eta (\vartheta ,0).$ Then A(0) and $A^*$
are adjoint operators. From Section 2, we know
that $\pm i \widetilde{\omega}\widetilde{\tau}_1$ are eigenvalues of A(0).
Thus, they are also eigenvalues of $A^*$. We first need to
compute the eigenvalues of A(0) and $A^*$ corresponding to
$i\widetilde{\omega}\widetilde{\tau}_1$ and $-i\widetilde{\omega}\widetilde{\tau}_1$,
respectively. Suppose that $q
(\vartheta) = (1,q_{1})^{T} e^{i \widetilde{\omega} \widetilde{\tau}_1 \vartheta}$ is
the eigenvector of A(0) corresponding to $i\widetilde{\omega}\widetilde{\tau}_1$, then  $A(0)q (\vartheta)$ = $ i\widetilde{\omega}\widetilde{\tau}_1
q (\vartheta).$ It follows from the definition of A(0)
and (4.2), (4.4), (4.5) that
$$\widetilde{\tau}_1\begin{pmatrix}
  _{i\widetilde{\omega}-\frac{\alpha^2(1-m)^2hx^*y^*}{\{1+\alpha (1-m)h x^*\}^2}+\frac{rx^*e^{-i\widetilde{\omega} \widetilde{\tau}_1}}{k}} & _{\frac{\alpha (1-m)x^*}{1+\alpha (1-m)hx^*}}\\
  _{-\frac{\theta \alpha (1-m)y^* e^{-i\widetilde{\omega} \widetilde{\tau}_2}}{\{1+\alpha (1-m)h x^*\}^2}} & _{i\widetilde{\omega}}
\end{pmatrix} q(0) = \begin{pmatrix}
  _{0} \\
  _{0}
\end{pmatrix}.$$ \noindent Thus, we can easily obtain\\ $$q(0) =
(1,q_{1})^{T},$$\\ where\\ $$q_{1} = \frac{\theta \alpha (1-m)y^* e^{-i\widetilde{\omega}\widetilde{\tau}_2}}{i\widetilde{\omega}\{1+\alpha (1-m)h x^*\}^2}.$$
Similarly, let $q^*(s) = D (1, {q_{1}}^*)^{T} $ is the eigenvector of $A^*$ corresponding to $- i
\widetilde{\omega} \widetilde{\tau}_1$. By the definition of $A^*$ and
(4.2),(4.3) and
(4.4), we can compute
$$q^*(s) = D (1, {q_{1}}^*)
         = D \left(1, -\frac{\alpha (1-m)x^*}{i\widetilde{\omega}\{1+\alpha (1-m)h x^*\}}\right).$$
In order to assure $\langle q^*(s), q(\vartheta) \rangle = 1$, we need
to determine the value of D. From (4.7), we have
\begin{eqnarray}\langle q^*(s) , q(q_{1}) \rangle &=& \overline{D} (1,
{\overline{q_{1}}}^*) (1,q_{1})^{T} - \int_{-\widetilde{\tau}_2}^{0} \int_{\xi =
0}^{\vartheta} \overline{D} (1,{\overline{q_{1}}}^*) e^{- i
\widetilde{\omega} (\xi - \vartheta)} d \eta (\vartheta)
(1,q_{1})^{T} e^{i \widetilde{\omega} \xi} d
\xi \nonumber\\
&=& \overline{D}\bigg\{\bigg(1-\frac{rx^*}{K}e^{-i\widetilde{\omega}\widetilde{\tau}_1}\bigg)+q_1 q_1^*\bigg(1+\frac{\theta \alpha (1-m)y^*e^{-i\widetilde{\omega} \widetilde{\tau}_2}}{\{1+\alpha(1-m)hx^*\}^2}\bigg)\bigg\}.\nonumber\end{eqnarray} \noindent Thus, we can choose D as
$$\overline{D} = \frac{1}{\bigg(1-\frac{rx^*}{K}e^{-i\widetilde{\omega}\widetilde{\tau}_1}\bigg)+q_1 q_1^*\bigg(1+\frac{\theta \alpha (1-m)y^*e^{-i\widetilde{\omega} \widetilde{\tau}_2}}{\{1+\alpha(1-m)hx^*\}^2}\bigg)}.$$ In the remainder
of the section, we use the same notations as in \cite{GH94}. We first
compute the coordinates to describe the center manifold $C_{0}$
at $\mu = 0$. Define
$$ z (t) = \langle q^*,x_{t}\rangle,~~~~~~~~W (t,\vartheta) = x_{t} (\vartheta) - 2
{\mbox Re} \{ z(t) q(\vartheta)\}.\eqno{(4.8)}$$  On the center
manifold $C_{0}$, we have
$$W (t,\vartheta) = W ( z(t),\overline{z} (t), \vartheta),$$
where
$$W (z,\overline{z}, \vartheta) = W_{20} (\vartheta) \frac{z^2}{2} + W_{11} (\vartheta) z \overline{z} + W_{02} (\vartheta) \frac{{\overline{z}}^2}{2} + W_{30} (\vartheta) \frac{z^3}{6} +......,\eqno{(4.9)}$$
z and $\overline{z}$ are local coordinates for center manifold
$C_{0}$ in the direction of $q^*$ and ${\overline{q}}^*.$ Note
that W is real if $x_{t}$ is real. We only consider real
solutions. For solution $x_{t} \in C_{0}$ of (4.6), since $\mu =
0,$ we have
$$ \dot{z} (t) = i \widetilde{\omega}\widetilde{\tau}_1z + {\overline{q}}^* (\vartheta) f (0,W (z,\overline{z},\vartheta) + 2 {\mbox Re}\{zq
(\vartheta)\}) ~~ ^{{\underline{\underline{\textmd{{def}}}}}} ~~i
\widetilde{\omega} \widetilde{\tau}_1 z + {\overline{q}}^* (0) f_{0}
(z,\overline{z}).$$ We rewrite this equation as
$$\dot{z} (t) = i \widetilde{\omega} \widetilde{\tau}_1z (t) + g (z,\overline{z}),$$
where
$$g (z,\overline{z}) = {\overline{q}}^* (0) f_{0} (z,\overline{z}) = g_{20} \frac{z^2}{2} + g_{11} z \overline{z} + g_{02} \frac{{\overline{z}}^2}{2} + g_{21} \frac{z^2 \overline{z}}{2} + .......\eqno{(4.10)}$$
It follows from (4.8) and (4.9) that
\begin{eqnarray}x_{t}(\vartheta) &=& W (t,\vartheta) + 2 {\mbox Re} \{z (t) q (t)\}\nonumber\\
&=& W_{20} (\vartheta) \frac{z^2}{2} + W_{11} (\vartheta) z
\overline{z} + W_{02} (\vartheta) \frac{{\overline{z}}^2}{2} +
(1,q_{1})^{T}
e^{i \widetilde{\omega} \widetilde{\tau}_1 \vartheta} z \nonumber\\
&+& (1,\overline{q_{1}})^{T} e^{- i \widetilde{\omega}\widetilde{\tau}_1
\vartheta} \overline{z} +..........
~~~~~~~~~~~~~~~~~~~~~~~~~~~~~~~~~~~~~~~~~~~~~~~(4.11)\nonumber
\end{eqnarray}
It follows together with (4.3) that
\begin{eqnarray}g(z,\overline{z}) &=&
{\overline{q}}^* (0) f (0,x_{t}) = \widetilde{\tau}_1\overline{D}
(1,{\overline{q_{1}}}^*) \begin{pmatrix}
  _{-\frac{r}{k}x_{1t}(0)x_{1t}(-\widetilde{\tau}_1)-\frac{\alpha (1-m)x_{1t}(0)x_{2t}(0)}{1+\alpha(1-m)hx_{1t}(0)}} \\
  _{\frac{\theta \alpha (1-m)x_{1t}(-\widetilde{\tau}_2)x_{2t}(0)}{1+\alpha(1-m)hx_{1t}(-\widetilde{\tau}_2)}}
\end{pmatrix} \nonumber\\
&=& -\frac{\widetilde{\tau}_1r\overline{D}}{k}[z^2e^{-i\widetilde{\omega}\widetilde{\tau}_1}+\overline{z}^2e^{i\widetilde{\omega}\widetilde{\tau}_1}+(e^{i\widetilde{\omega}\widetilde{\tau}_1}+
e^{-i\widetilde{\omega}\widetilde{\tau}_1})z\overline{z}+\{2(W_{11}^{1}(-\widetilde{\tau}_1)+W_{11}^{1}(0)e^{-i\widetilde{\omega}
\widetilde{\tau}_1})\nonumber\\
&+&W_{20}^{1}(-\widetilde{\tau}_1)+W_{20}^{1}(0)e^{i\widetilde{\omega}\widetilde{\tau}_1}\}\frac{z^2\overline{z}}{2}]
-\widetilde{\tau}_1\overline{D}\alpha (1-m)[q_1z^2+\overline{q_1}\overline{z}^2+(q_1+\overline{q_1})z\overline{z}\nonumber\\
&+&\frac{z^2\overline{z}}{2}\{2(q_1W_{11}(0)+W_{11}^2(0))
+(\overline{q_1}W_{20}^1(0)+W_{20}^2(0))\}]\nonumber\\
&+&2\widetilde{\tau}_1\overline{D}\alpha^2(1-m)^2h(\overline{q_1}+2q_1)\frac{z^2\overline{z}}{2}+\widetilde{\tau}_1\overline{D}\overline{q_1^*}\theta \alpha (1-m)[q_1e^{-i\widetilde{\omega}\widetilde{\tau}_2}z^2+\overline{q_1}e^{i\widetilde{\omega}\widetilde{\tau}_2}\overline{z}^2\nonumber\\
&+&(\overline{q_1}e^{-i\widetilde{\omega}\widetilde{\tau}_2}+q_1e^{i\widetilde{\omega}\widetilde{\tau}_2})z\overline{z}+\{2(W_{11}^2(0)e^{-i\widetilde{\omega}\widetilde{\tau}_2}+W_{11}^1(-\widetilde{\tau}_2)q_1)\nonumber\\
&+&(W_{20}^2)(0)e^{i\widetilde{\omega}\widetilde{\tau}_2}
+\overline{q_1}W_{20}^1(-\widetilde{\tau}_2)\}\frac{z^2\overline{z}}{2}]\nonumber\\
&-&2\widetilde{\tau}_1\overline{D}\overline{q_1^*}\theta \alpha^2 (1-m)^2 h (\overline{q_1}e^{-2i\widetilde{\omega}\widetilde{\tau}_2}+2q_1)\frac{z^2\overline{z}}{2}........~~~~~~~~~~~~~~~~~~~~~~~~~~~~~~~~~~ (4.12)\nonumber
\end{eqnarray}
Comparing the coefficients with $(4.10)$, we have
\begin{eqnarray}
g_{20} &=& 2\widetilde{\tau}_1\overline{D}[-\frac{r}{k}e^{-i\widetilde{\omega}\widetilde{\tau}_1}+\alpha(1-m)q_1(\theta\overline{q_1}^*e^{-i\widetilde{\omega}\widetilde{\tau}_2}-1)],\nonumber\\
g_{11} &=& 2 \widetilde{\tau}_1\overline{D} [-\frac{r}{k}(e^{i\widetilde{\omega}\widetilde{\tau}_1}+e^{-i\widetilde{\omega}\widetilde{\tau}_1})+\alpha(1-m)\{\theta\overline{q_1}^*(\overline{q_1}
e^{-i\widetilde{\omega}\widetilde{\tau}_2}+q_1e^{i\widetilde{\omega}\widetilde{\tau}_2})-Re\{q_1\}\}],\nonumber\\
g_{02} &=& 2\widetilde{\tau}_1\overline{D}[-\frac{r}{k}e^{i\widetilde{\omega}\widetilde{\tau}_1}+\alpha(1-m)\overline{q_1}(\theta\overline{q_1}^*
e^{i\widetilde{\omega}\widetilde{\tau}_2}-1)],\nonumber\\
g_{21} &=& -\frac{\widetilde{\tau}_1r\overline{D}}{k}\{2(W_{11}^1(-\widetilde{\tau}_1)+W_{11}^2(0)e^{-i\widetilde{\omega}
\widetilde{\tau}_1})+W_{20}^1(0)(-\widetilde{\tau}_1)+W_{20}^2(0)e^{i\widetilde{\omega}\widetilde{\tau}_1}\}\nonumber\\ &-&\widetilde{\tau}_1\overline{D}\alpha(1-m)\{2(q_1W_{11}^1(0)+W_{11}^2(0))+(\overline{q_1}W_{20}^1(0)+W_{20}^2(0))\}\nonumber\\
&+& 2\widetilde{\tau}_1\overline{D}\alpha^2(1-m)^2h(\overline{q_1}+2q_1)
+\widetilde{\tau}_1\overline{D}\overline{q_1}^*\theta\alpha(1-m)\{W_{11}^2(0)e^{-i\widetilde{\omega}\widetilde{\tau}_2}+W_{11}^1
(-\widetilde{\tau}_2q_1)\nonumber\\
&+&(W_{20}^2(0)e^{i\widetilde{\omega}\widetilde{\tau}_2}+\overline{q_1}W_{20}^1(-\widetilde{\tau}_2))\}-2\widetilde{\tau}_1\overline{D}
\overline{q_1}^*
\theta\alpha^2(1-m)^2h(\overline{q_1}e^{-2i\widetilde{\omega}\widetilde{\tau}_2}+2q_1)
.~~~~~~~(4.13)\nonumber
\end{eqnarray}
Since there are $W_{20} (\vartheta)$ and $W_{11} (\vartheta)$ in
$g_{21}$, we still need to compute them.

\noindent From (4.3) and (4.8), we have
\begin{eqnarray}\dot{W} = \dot{x_{t}} - \dot{z} - \dot{\overline{z}}
\overline{q} &=& \left\{\begin{array}{ll}
  A W - 2 {\mbox Re} \{{\overline{q}}^* (0) f_{0 q} (\vartheta)\}, & \vartheta \in [-1,0), \\
  A W - 2 {\mbox Re} \{{\overline{q}}^* (0) f_{0 q}
  (\vartheta)\} + f_{0}, & \vartheta = 0,
\end{array}\right.\nonumber\\
 & ^{\underline{\underline{\textmd{def}}}}& A W + H (z, \overline{z},
 \vartheta),~~~~~~~~~~~~~~~~~~~~~~~~~~~~~~~~~~~~~~~~~~~~~~~~~~~~~~~(4.14)\nonumber\end{eqnarray}
where
$$H (z, \overline{z}, \vartheta) = H_{20} (\vartheta) \frac{z^2}{2} + H_{11} (\vartheta) z \overline{z} + H_{02} (\vartheta) \frac{{\overline{z}}^2}{2} + .....\eqno{(4.15)}$$
Substituting the corresponding series into (4.14) and comparing
the coefficients, we obtain
$$(A - 2 i \widetilde{\omega}\widetilde{\tau}_1 ) W_{20} (\vartheta) = - H_{20}(\vartheta),~~~~~~~A W_{11} (\vartheta) = - H_{11} (\vartheta)......\eqno{(4.16)}$$
From (4.14), we know that for $\vartheta \in [-1,0),$
$$H (z, \overline{z}, \vartheta) = - {\overline{q}}^* (0) f_{0} q (\vartheta) - {q}^* (0)\overline{f}_{0} (\vartheta) = - g (z,\overline{z}) q (\vartheta) - \overline{g} (z,\overline{z}) \overline{q}(\vartheta),\eqno{(4.17)}.$$
Comparing the coefficients with (4.15), we get
$$H_{20} (\vartheta) = - g_{20} q (\vartheta) - \overline{g}_{02} \overline{q
} (\vartheta)$$ and
$$H_{11} (\vartheta) = - g_{11} q (\vartheta) - \overline{g}_{11} \overline{q
} (\vartheta).$$ From (4.16) and (4.18) and the definition of A, it
follows that
$$\dot{W}_{20} (\vartheta) = 2 i \widetilde{\omega} \widetilde{\tau}_1W_{20} (\vartheta) + g_{20} q (\vartheta) + \overline{g}_{02} \overline{q} (\vartheta) .$$
Notice that $q (\vartheta) = (1,x)^{T} e^{i \widetilde{\omega}\widetilde{\tau}_1
\vartheta} ,$ hence
$$W_{20} (\vartheta) = \frac{i g_{20}}{\widetilde{\omega} } q(0) e^{i \widetilde{\omega} \widetilde{\tau}_1\vartheta} + \frac{i \overline{g}_{02}}{3 \widetilde{\omega} } \overline{q} (0) e^{- i \widetilde{\omega}\widetilde{\tau}_1 \vartheta} + E_{1}^{'} e^{2 i \widetilde{\omega}\widetilde{\tau}_1 \vartheta},\eqno{(4.20)}$$
where $E_{1}^{'} = (E_{1}^{(1)},E_{1}^{(2)}) \in R^{2}$ is a
constant vector. \noindent Similarly, from (4.16) and (4.19), we
obtain
$$W_{11} (\vartheta) = - \frac{i g_{11}}{\widetilde{\omega} } q(0) e^{i \widetilde{\omega}\widetilde{\tau}_1  \vartheta} + \frac{i \overline{g}_{11}}{\widetilde{\omega} } \overline{q} (0) e^{- i \widetilde{\omega} \widetilde{\tau}_1 \vartheta} + E_{2}^{'},\eqno{(4.21)}$$
where $E_{2}^{'} = (E_{2}^{(1)},E_{2}^{(2)}) \in R^{2}$ is a
constant vector. \noindent In what follows, we shall seek
appropriate $E_{1}$ and $E_{2}$. From the definition of A and
(3.16), we obtain
$$\int_{-1}^{0} d \eta (\vartheta) W_{20} (\vartheta) = 2 i \widetilde{\omega}\widetilde{\tau}_1 W_{20} (0) - H_{20}
(0)\eqno{(4.22)}$$ and
$$\int_{-1}^{0} d \eta (\vartheta) W_{11} (\vartheta) = - H_{11} (0),\eqno{(4.23)}$$
where $\eta(\vartheta) = \eta(0,\vartheta).$ By (4.14), we have
$$H_{20} (0) = - g_{20} q(0) - \overline{g}_{02} \overline{q} (0) +2\widetilde{\tau}_1 \begin{pmatrix}
  _{-\bigg(\frac{r}{k}e^{-i\widetilde{\omega}\widetilde{\tau}_1}+\alpha(1-m)q_1\bigg)} \\
  _{\theta \alpha (1-m)e^{-i\widetilde{\omega}\widetilde{\tau}_2}q_1}
\end{pmatrix}\eqno{(4.24)}$$ and
$$H_{11} (0) = - g_{11} q(0) - \overline{g}_{11} \overline{q} (0) + 2\widetilde{\tau}_1 \begin{pmatrix}
  _{-\frac{r}{k}(e^{-i\widetilde{\omega}\widetilde{\tau}_1}+e^{i\widetilde{\omega}\widetilde{\tau}_1})-\alpha(1-m)Re\{q_1\}} \\
  _{\alpha(1-m)(\overline{q_1}e^{-i\widetilde{\omega}\widetilde{\tau}_2}+q_1e^{i\widetilde{\omega}\widetilde{\tau}_2})}
\end{pmatrix}.\eqno{(4.25)}$$ Substituting (4.20) and (4.24) into (4.22) and
noticing that
$$\left(i \widetilde{\omega}\widetilde{\tau}_1 I - \int_{-1}^{0} e^{i \widetilde{\omega}\widetilde{\tau}_1\vartheta} d \eta (\vartheta)\right) q (0) = 0$$
and
$$\left(- i \widetilde{\omega} \widetilde{\tau}_1I - \int_{-1}^{0} e^{- i \widetilde{\omega} \widetilde{\tau}_1 \vartheta} d \eta (\vartheta)\right) \overline{q} (0) = 0,$$
we obtain
$$\left(2 i \widetilde{\omega}\widetilde{\tau}_1 I - \int_{-1}^{0} e^{2 i \widetilde{\omega}\widetilde{\tau}_1 \vartheta} d \eta (\vartheta)\right) E_{1}^{'} = 2\widetilde{\tau}_1 \begin{pmatrix}
  _{-\bigg(\frac{r}{k}e^{-i\widetilde{\omega}\widetilde{\tau}_1}+\alpha(1-m)q_1\bigg)} \\
  _{\theta \alpha (1-m)e^{-i\widetilde{\omega}\widetilde{\tau}_2}q_1}
\end{pmatrix}. $$ This leads to
$$\begin{pmatrix}
  _{2i\widetilde{\omega}-\frac{\alpha^2(1-m)^2hx^*y^*}{\{1+\alpha (1-m)h x^*\}^2}+\frac{rx^*}{k}}e^{-i\widetilde{\omega} \widetilde{\tau}_1} & _{\frac{\alpha (1-m)x^*}{1+\alpha (1-m)hx^*}}\\
  _{-\frac{\theta \alpha (1-m)y^* e^{-i\widetilde{\omega} \widetilde{\tau}_2}}{\{1+\alpha (1-m)h x^*\}^2}} & _{2i\widetilde{\omega}}
\end{pmatrix} E_{1}^{'} = 2 \begin{pmatrix}
  _{-\bigg(\frac{r}{k}e^{-i\widetilde{\omega}\widetilde{\tau}_1}+\alpha(1-m)q_1\bigg)} \\
  _{\theta \alpha (1-m)e^{-i\widetilde{\omega}\widetilde{\tau}_2}q_1}
\end{pmatrix}.$$
Therefore, it follows that
$$E_{1}^{(1)} = \frac{2}{A}~ \vline \begin{array}{cc}
  _{-\bigg(\frac{r}{k}e^{-i\widetilde{\omega}\widetilde{\tau}_1}+\alpha(1-m)q_1\bigg)} & _{\frac{\alpha (1-m)x^*}{1+\alpha (1-m)hx^*}} \\
  _{\theta \alpha (1-m)e^{-i\widetilde{\omega}\widetilde{\tau}_2}q_1} & _{2i\widetilde{\omega}}
\end{array}\vline $$\\
and\\
$$E_{1}^{(2)} = \frac{2}{A}~ \vline \begin{array}{cc}
  _{2i\widetilde{\omega}-\frac{\alpha^2(1-m)^2hx^*y^*}{\{1+\alpha (1-m)h x^*\}^2}+\frac{rx^*}{K}}e^{-i\widetilde{\omega} \widetilde{\tau}_1} & {-\bigg(\frac{r}{k}e^{-i\widetilde{\omega}\widetilde{\tau}_1}+\alpha(1-m)q_1\bigg)}\\
 _{-\frac{\theta \alpha (1-m)y^* e^{-i\widetilde{\omega} \widetilde{\tau}_2}}{\{1+\alpha (1-m)h x^*\}^2}} & {\theta \alpha (1-m)e^{-i\widetilde{\omega}\widetilde{\tau}_2}q_1}
\end{array} \vline~.$$
where
$$A = ~\vline \begin{array}{cc}
  _{2i\widetilde{\omega}-\frac{\alpha^2(1-m)^2hx^*y^*}{\{1+\alpha (1-m)h x^*\}^2}+\frac{rx^*}{k}}e^{-i\widetilde{\omega} \widetilde{\tau}_1} & _{\frac{\alpha (1-m)x^*}{1+\alpha (1-m)hx^*}}\\
  _{-\frac{\theta \alpha (1-m)y^* e^{-i\widetilde{\omega} \widetilde{\tau}_2}}{\{1+\alpha (1-m)h x^*\}^2}} & _{2i\widetilde{\omega}}
\end{array}\vline~.$$
Similarly, substituting (4.21) and (4.25) into (4.23), we get
$$\begin{pmatrix}
  _{-\frac{\alpha^2(1-m)^2hx^*y^*}{\{1+\alpha (1-m)h x^*\}^2}+\frac{rx^*}{k}} & _{\frac{\alpha (1-m)x^*}{1+\alpha (1-m)hx^*}}\\
  _{-\frac{\theta \alpha (1-m)y^* }{\{1+\alpha (1-m)h x^*\}^2}} & _{0}
\end{pmatrix} E_{2}^{'} = 2 \begin{pmatrix}
 _{-\frac{r}{k}(e^{-i\widetilde{\omega}\widetilde{\tau}_1}+e^{i\widetilde{\omega}\widetilde{\tau}_1})-\alpha(1-m)Re\{q_1\}} \\
  _{\alpha(1-m)(\overline{q_1}e^{-i\widetilde{\omega}\widetilde{\tau}_2}+q_1e^{i\widetilde{\omega}\widetilde{\tau}_2})}
\end{pmatrix}.$$
It follows that
$$E_{2}^{(1)} = \frac{2}{B^{'}}~ \vline \begin{array}{cc}
 -\frac{r}{k}(e^{-i\widetilde{\omega}\widetilde{\tau}_1}+e^{i\widetilde{\omega}\widetilde{\tau}_1})-\alpha(1-m)Re\{q_1\} & \frac{\alpha (1-m)x^*}{1+\alpha (1-m)hx^*}\\
  \alpha(1-m)(\overline{q_1}e^{-i\widetilde{\omega}\widetilde{\tau}_2}+q_1e^{i\widetilde{\omega}\widetilde{\tau}_2}) & 0
\end{array} \vline $$
and
$$E_{2}^{(2)} = \frac{2}{B^{'}}~ \vline \begin{array}{cc}
 -\frac{\alpha^2(1-m)^2hx^*y^*}{\{1+\alpha (1-m)h x^*\}^2}+\frac{rx^*}{K} & -\frac{r}{k}(e^{-i\widetilde{\omega}\widetilde{\tau}_1}+e^{i\widetilde{\omega}\widetilde{\tau}_1})-\alpha(1-m)Re\{q_1\} \\
-\frac{\theta \alpha (1-m)y^* }{\{1+\alpha (1-m)h x^*\}^2} & \alpha(1-m)(\overline{q_1}e^{-i\widetilde{\omega}\widetilde{\tau}_2}+q_1e^{i\widetilde{\omega}\widetilde{\tau}_2})
\end{array} \vline~,$$
where
$$B^{'} =~ \vline \begin{array}{cc}
   _{-\frac{\alpha^2(1-m)^2hx^*y^*}{\{1+\alpha (1-m)h x^*\}^2}+\frac{rx^*}{k}} & _{\frac{\alpha (1-m)x^*}{1+\alpha (1-m)hx^*}}\\
  _{-\frac{\theta \alpha (1-m)y^* }{\{1+\alpha (1-m)h x^*\}^2}} & _{0}
\end{array}\vline~.$$
Thus, we can determine $W_{20} (\vartheta)$ and $W_{11} (\vartheta)$
from (4.20) and (4.21). Furthermore, $g_{21}$ in (4.13) can be
expressed by the parameters and delay. Thus, we can compute the
following values:
\begin{eqnarray}
&&c_{1} (0) = \frac{i}{2 \widetilde{\omega} \tau_{1_{0}}} (g_{20} g_{11} - 2 |g_{11}|^2 - \frac{|g_{02}|^2}{3}) + \frac{g_{21}}{2},\nonumber\\
&&\mu_{2} = - \frac{{\mbox Re} \{c_{1} (0)\}}{{\mbox Re} \{\lambda^{'} (\tau_{1_{0}})\}}, \nonumber\\
&&\beta_{2} = 2 {\mbox Re} (c_{1} (0)), \nonumber\\
&&T_{2} = - \frac{{\mbox Im} \{c_{1} (0)\} + \mu_{2} {\mbox Im}
\{\lambda^{'} (\tau_{1_{0}})\}}{\widetilde{\omega}
  \tau_{1_{0}}}.\nonumber
\end{eqnarray}
which determine the qualities of bifurcating periodic solution in
the center manifold at the critical value $\tau_{1_{0}}$.

\noindent {\bf Lemma 3.1} $\mu_{2}$ determines the direction of the Hopf bifurcation. If
$\mu_{2} > 0  ~(\mu_{2} < 0)$ then the Hopf bifurcation is
supercritical (subcritical) and the bifurcating periodic solutions
exist for $\tau > \tau_{1_{0}} ~ (\tau < \tau_{1_{0}}). ~\beta_{2}$ determines the stability of the bifurcating periodic
solutions: the bifurcating periodic solutions are stable (unstable)
if $\beta_{2} < 0 ~ (\beta_{2} > 0)$. $T_{2}$ determines the
period of the bifurcating periodic solutions: the period increases
(decreases) if $T_{2} > 0 ~ (< 0).$

\section{Simulation results}
In this section, we give some numerical simulations to illustrate
the analytical results observed in the previous sections. For
illustration purpose, we consider the parameter values $r=2.65,~k=898,~\alpha=0.045,~m=0.45,~h=0.0437,~\theta=0.215,~d=1.06$ with initial value $(30, 5.83)$.\\

\begin{figure}
\centering
\includegraphics[width=5in,height=3in]{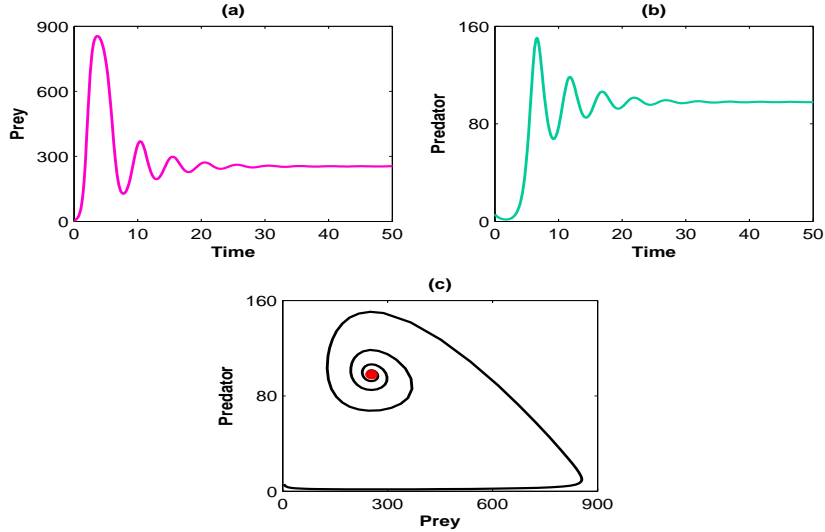}
\caption{The coexistence
equilibrium point $E^{*}=(253.9056,97.8867)$ is locally asymptotically stable when $\tau_1=0=\tau_2$. Figs. (a) and (b) are the time evolutions of the model system (1.3) and Fig. (c) is the corresponding phase plane.}
\label{fig2}
\end{figure}

\begin{figure}
\centering
\includegraphics[width=5in,height=4in]{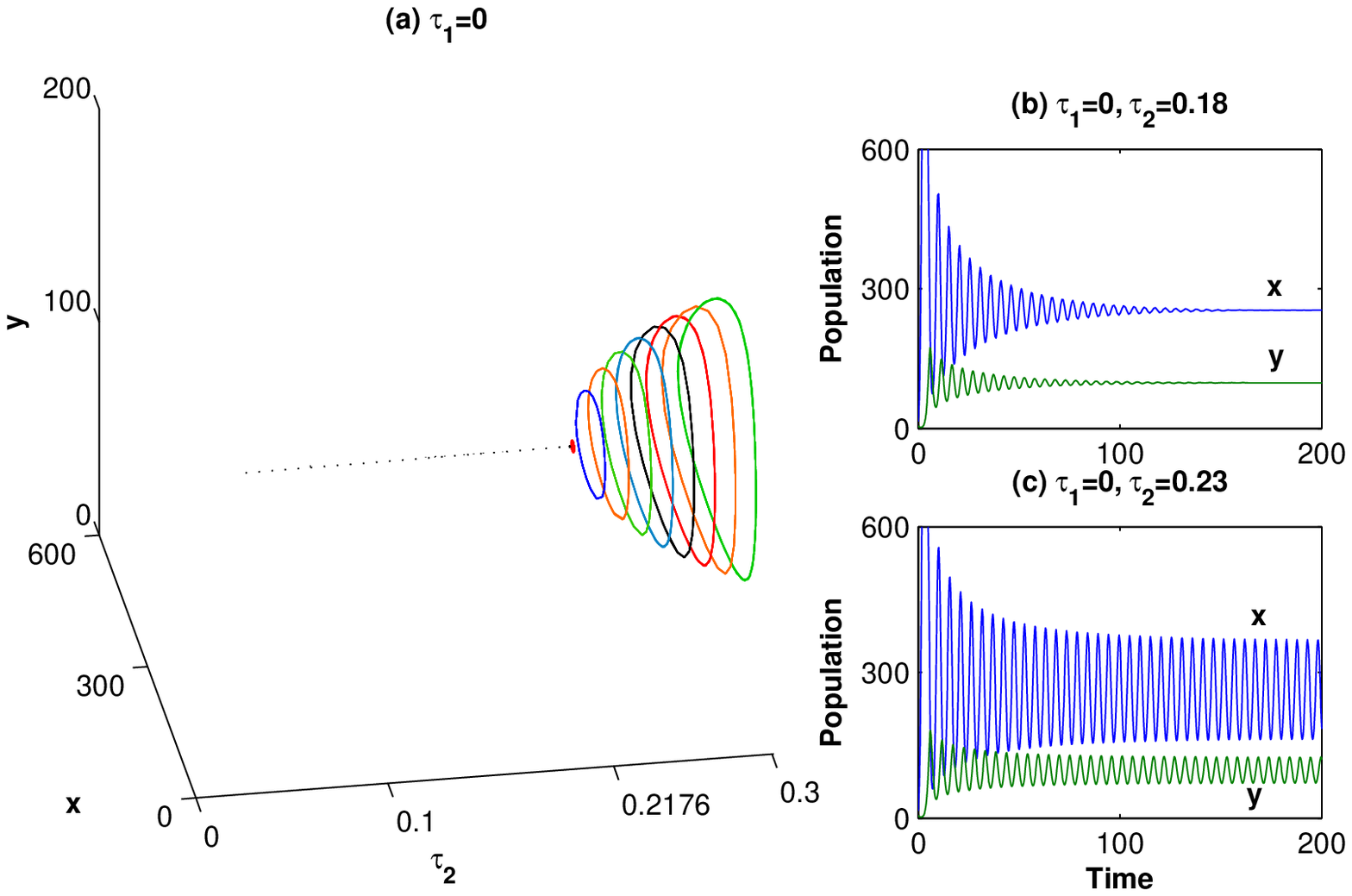}
\caption{Bifurcation diagram (Figure 2(a))of system (1.3) with respect to the
bifurcation pa\-ra\-me\-ter $\tau_2$ is drawn in the three-dimensional space $(\tau_2, x, y)$ when $\tau_1=0$.  Figure 2(b) shows the time evolution of the system when $\tau_2 (=0.18)$ is less than its critical value $0.2176$. Figure 2(c) indicates the time evolution of the system when $\tau_2 (=0.23)$ is greater than its critical value $0.2176$. These figures show that the coexistence equilibrium is stable for $\tau_2 <0.2176$, unstable for $\tau_2 >0.2176$ and Hopf-bifurcation occurs at $\tau=\tau_{2_0}=0.2176$ with $\tau_1=0$.}
\label{fig2}
\end{figure}

\begin{figure}
\centering
\includegraphics[width=5in,height=3in]{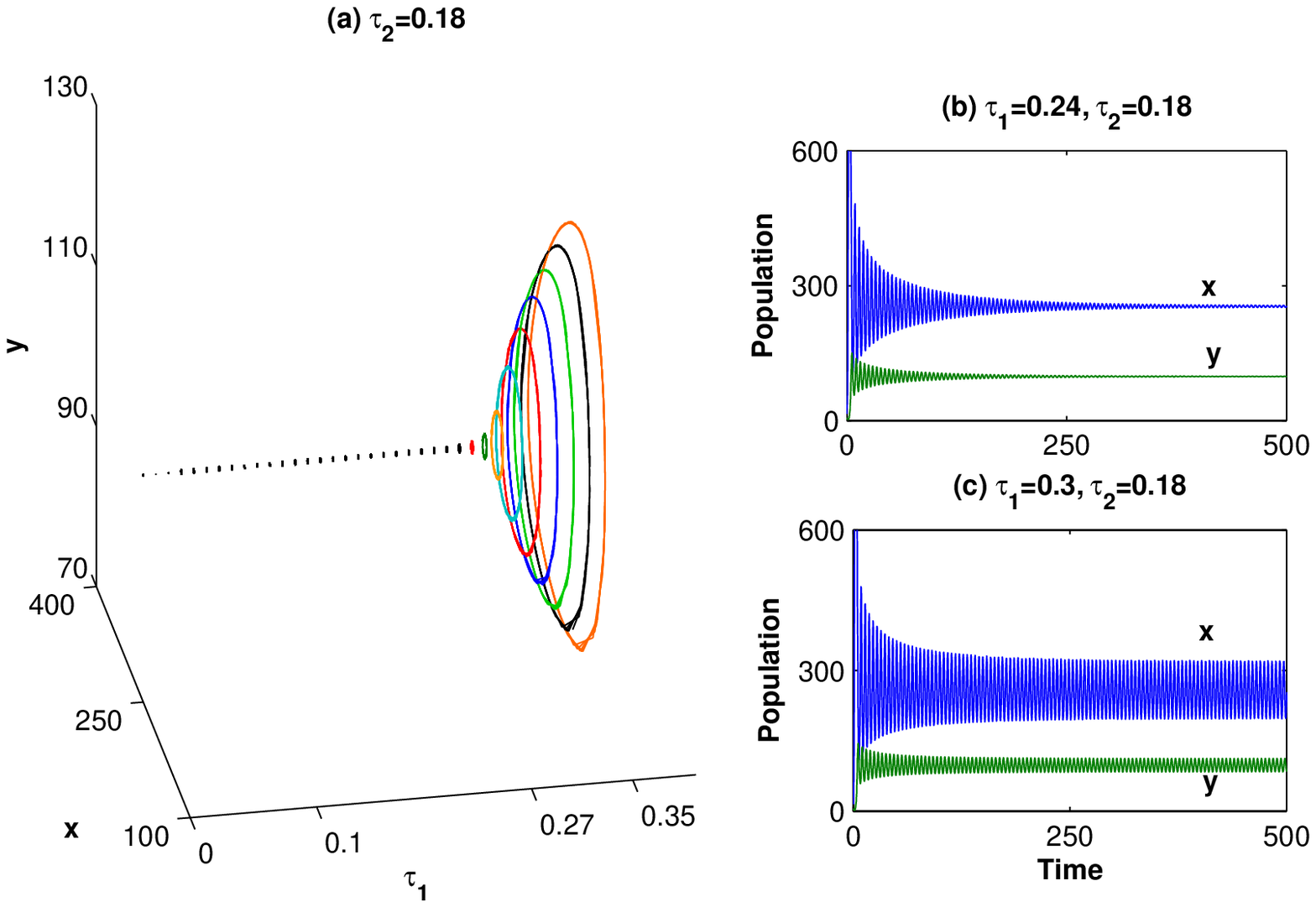}
\caption{Bifurcation diagram (Figure 3(a))of system (1.3) with respect to the
bifurcation pa\-ra\-me\-ter $\tau_2$ is drawn in the three-dimensional space $(\tau_1, x, y)$ when $\tau_2=0.18$.  Figure 3(b) shows the time evolution of the system when $\tau_1 (=0.24)$ is less than its critical value $0.27$. Figure 3(c) indicates the time evolution of the system when $\tau_1 (=0.3)$ is greater than its critical value $0.27$. These figures show that the coexistence equilibrium is stable for $\tau_1 <0.27$, unstable for $\tau_1 >0.27$ and Hopf-bifurcation occurs at $\tau=\tau_{1_0}=0.27$ with $\tau_2=0.18$.}
\label{fig2}
\end{figure}

%\begin{figure}
%\centering
%\includegraphics[width=5in,height=3in]{./fig4.eps}
%\caption{Stable and unstable regions of the system $(1.3)$ in the $\tau_2-c$ in the Case II. Here $\tau_1=0$ and other
%parameters are as in the Table 1.}
%\label{fig2}
%\end{figure}

%\begin{figure}
%\centering
%\includegraphics[width=5in,height=4in]{./fig5.eps}
%\caption{Figures in the first panel (a-c) depict that the system $(1.3)$ is stable for
%$\tau_{1}=0.24 ~(<\tau_{1_{0}}=0.27)$ and figures in the second panel (d-f) depict that the system  is unstable for
%$\tau_{1}=0.3 ~(>\tau_{1_{0}}=0.27)$ taking $\tau_{2}=0.18\in[0, \tau_{2_{0}}]$ in each case.
%Here $\tau_2=0.18$ and other parameters are as in the Table 1.}
%\label{fig2}
%\end{figure}

\begin{figure}
\centering
\includegraphics[width=5in,height=4in]{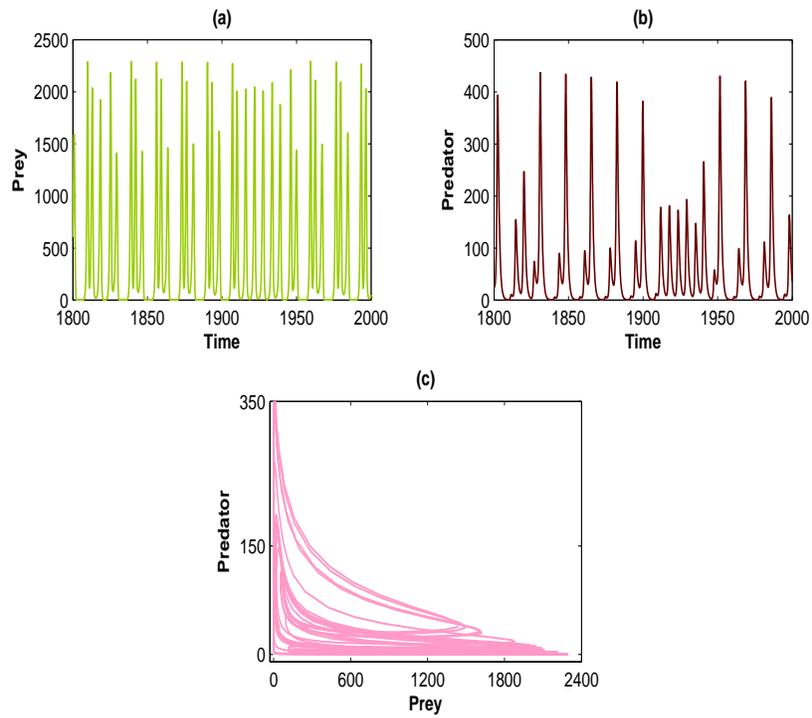}
\caption{Irregular periodic behavior of the system $(1.3)$ when $\tau_1=0.7 ~(>\tau_{1_{0}})$
and $\tau_2=0.8 ~(>\tau_{2_{0}})$. Figs. (a) and (b) are the time series of the prey and predator populations and Fig. (c) is the corresponding phase
diagram.}
\label{fig2}
\end{figure}

The above parameter set satisfies all the conditions in Lemma 2.1 and
consequently the system equation $(1.3)$ converges to the coexistence
equilibrium point $E^{*}(253.9056,97.8867)$ when $\tau_{1}=0$,
$\tau_{2}=0$ (Fig. 1).
For the Case II, where $\tau_{1}=0$ and $\tau_{2}\neq 0$, we compute $\omega_{0}=1.2345$ and
$\tau_{2_{0}}=0.2176$. Therefore, by Lemma $2.2$,
$E^{*}$ is asymptotically stable for
$\tau_{2}<\tau_{2_{0}}$ and unstable for
$\tau_{2}>\tau_{2_{0}}$. The system experiences a Hopf bifurcation around $E^{*}$
when $\tau_{2}=\tau_{2_{0}}$. This behavior is depicted in the Fig. 2. Figure 2 (b-c) depict,
respectively, that the system $(1.3)$ is stable for $\tau_{2}=0.18 ~(<\tau_{2_{0}}=0.2176)$ and unstable for
$\tau_{2}=0.23 ~(>\tau_{2_{0}}=0.2176)$ with $\tau_{1}=0$. The bifurcation diagram (Fig. 2(a)) clearly
demonstrates the system behavior for different values of $\tau_{2}$. It shows that
the system remains stable if $\tau_{2}<\tau_{2_{0}}$; but
the instability sets in through periodic oscillations if $\tau_{2}>\tau_{2_{0}}$. The system undergoes a Hopf bifurcation when $\tau_{2}=\tau_{2_{0}}=0.2176$. %The stability region is represented in Fig. 4 for wide ranges of $\tau_2$ and $c$.\\
\begin{figure}
\centering
\includegraphics[width=5in,height=4in]{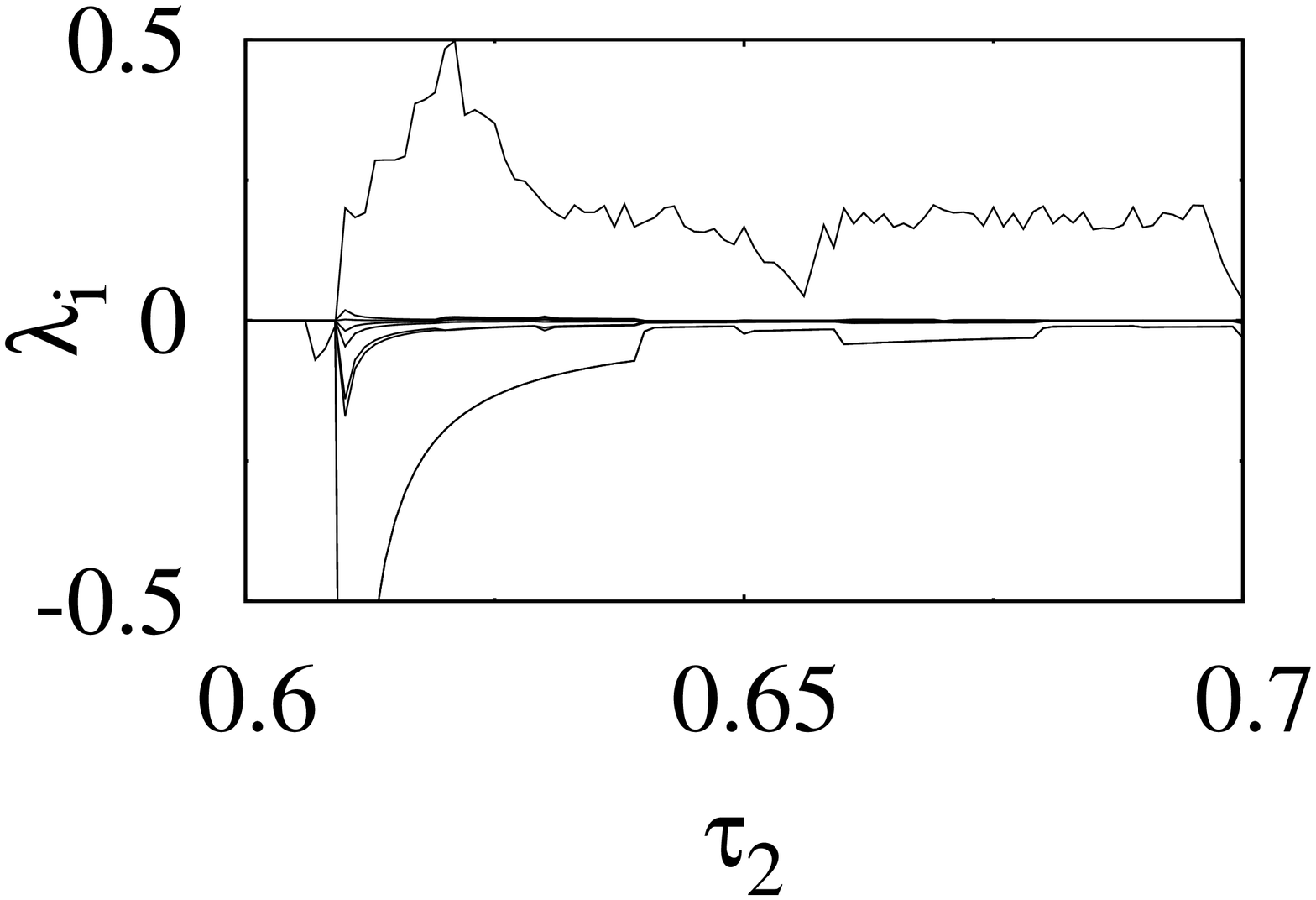}
\caption{Lyapunov exponents test of the chaotic solutions of the system (1.3).
Parameters are as in the Fig. 4. Here the largest few exponents have been plotted as a function of $\tau_2$. The highest one, namely $\lambda_1$ is positive for a large range of $\tau_2$ indicating chaos.}
\label{fig2}
\end{figure}
For the Case III,  where $\tau_{1}\neq 0$ and $\tau_{2}\neq 0$, we vary $\tau_{1}$ keeping the value of $\tau_{2}$ within its stability range $(0,0.2176)$. Choosing $\tau_{2}=0.18$, we obtain $\omega_{1}=1.1095$ as a root of the equation $(2.7)$. In this case, the value of $\tau_{1_{0}}$ becomes $0.27$ and the system $(1.3)$ (following the Lemma 2.3) is locally stable (unstable) whenever $\tau_1<0.27$ ($>0.27)$. A Hopf bifurcation occurs when $\tau_1=0.27$. Fig. 3(b) shows that the system (1.3) is locally asymptotically stable for $\tau_{1}=0.24<0.27$; and Fig. 3(c) show that the system (1.3) is unstable for $\tau_{1}=0.3>0.27$. The bifurcation diagram (Fig. 3(a)) clearly
demonstrates the system behavior for different values of $\tau_{1}$.
\begin{figure}
\centering
\includegraphics[width=5in,height=4in]{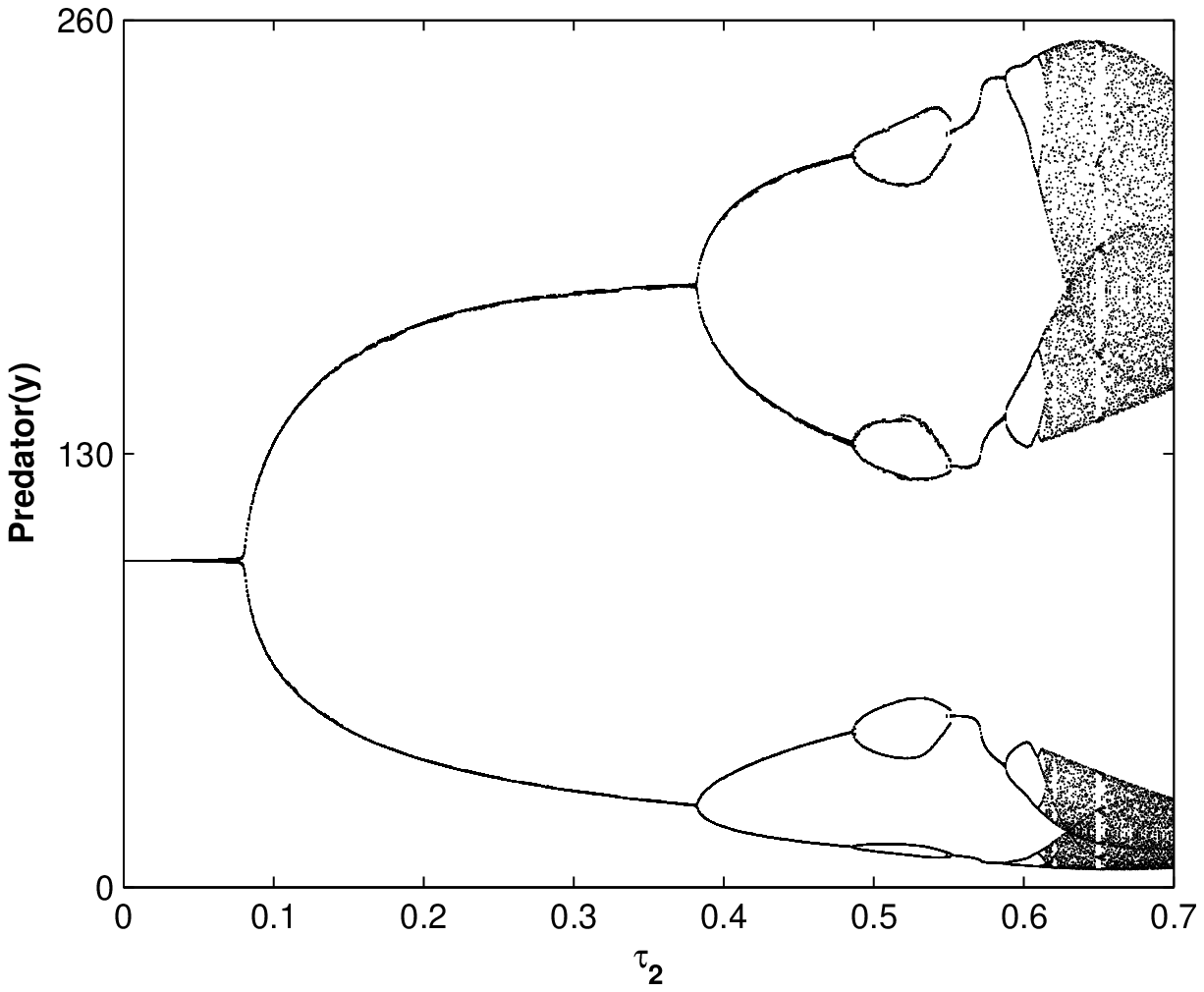}
\caption{Bifurcation diagram of the predator population when the delay parameter $\tau_2$ is smoothly varied with fixed $\tau_1=0.5$. Other parameters are as in Figure 1.}
\label{fig2}
\end{figure}
In Case IV, when $\tau_{2}=0, \tau_{1}\neq 0$, one can compute $\bar{\omega}_{0}=1.6095 $ and the corresponding critical value
of $\tau_1$ as
$\tau_{1_{0}}=0.6167$. Therefore, the coexistence equilibrium
$E^*(x^{*},y^{*})$ is locally asymptotically stable for $\tau_{1}<\bar{\tau}_{1_{0}}=0.6167$ and it is unstable for
$\tau_{1}>\bar{\tau}_{1_{0}}=0.6167$. Similarly for Case V, we take any value of $\tau_{1}$ from its stability range $[0, 0.6167)$, say $\tau_{1}=0.45$, and consider $\tau_{2}$ as a free parameter. One finds that $\bar{\omega}_1=1.6468$ and the corresponding critical value of $\tau_2$ is $\bar{\tau}_{2_0}=0.091$. Thus, for any fixed stable value of $\tau_{1}$, the system exhibits stable behavior around $E^*(x^{*},y^{*})$ for $\tau_{2}<\bar{\tau}_{2_{0}}$ and unstable oscillatory behavior for $\tau_{2}>\bar{\tau}_{2_{0}}$. \\
%\begin{figure}
%\centering
%\includegraphics[width=4in,height=3in]{./fig9.eps}
%\caption{Phase plane of the system $(1.3)$ with
%different $\tau_2 ~(>\tau_{2_{0}})$ but fixed $\tau_1 ~(>\tau_{1_{0}})$.
%(a) T period:
%$\tau_1=0.5, \tau_2=0.3$,
%(b) 2T period: $\tau_1=0.5, \tau_2=0.4$ and (c) 4T period: $\tau_1=0.5, \tau_2=0.53$. Other parameters
%are as in Table 1.}
%\label{fig2}
%\end{figure}

\begin{figure}
\centering
\includegraphics[width=4in,height=3in]{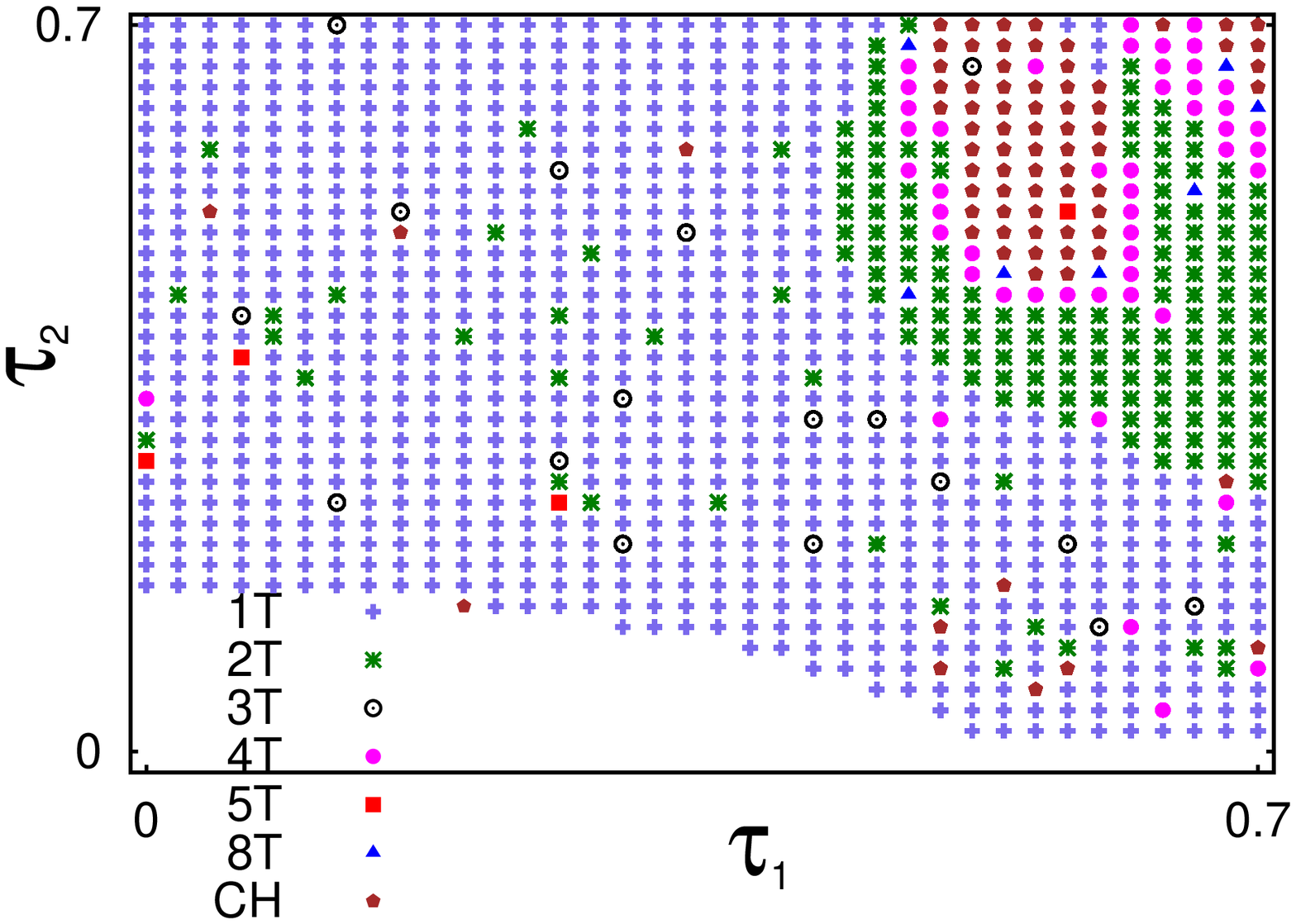}
\caption{Stability regions of the system $(1.3)$ are depicted in the $\tau_1\tau_2$-plane for  fixed $m ~(=0.45)$. Different dynamical features are depicted by different indices which are given in the figure.
Remaining parameters are as in the Figure 1.}
\label{fig2}
\end{figure}

Using Lemma 3.1, one can determine the values of $c_{1}(0),\mu_2,~ \beta_2$ and
$T_2$ as $c_{1}(0)=-0.000069936-0.00046189i,$ $\mu_2=0.000033112 (>0)$, $\beta_2=-0.000013987$ $(<0)$ and
$T_2=0.000023056 (>0)$. Since $\mu_2>0$ and $\beta_2 <0$, the Hopf bifurcation is supercritical and stable. Also, the period of the bifurcating periodic
solutions increases with $\tau_1$ (as $T_{2} > 0$), where $\tau_2$ is kept fixed in its stable region.\\

So far we have observed the behavior of the system when the delay parameters are within or slightly above the critical values. One interesting topic in the delay-induced system is to study the dynamical behavior of the system when the delay parameters are far away from their critical values, or they assume large values.
%For example, what type of dynamical behavior the system (1.3) will exhibit when the delay parameters assume values significantly larger than their respective critical values.
To observe the dynamics, we have simulated our model system (1.3) for larger values of $\tau_1$ and $\tau_2$.\\
\begin{figure}
\centering
\includegraphics[width=5in,height=4in]{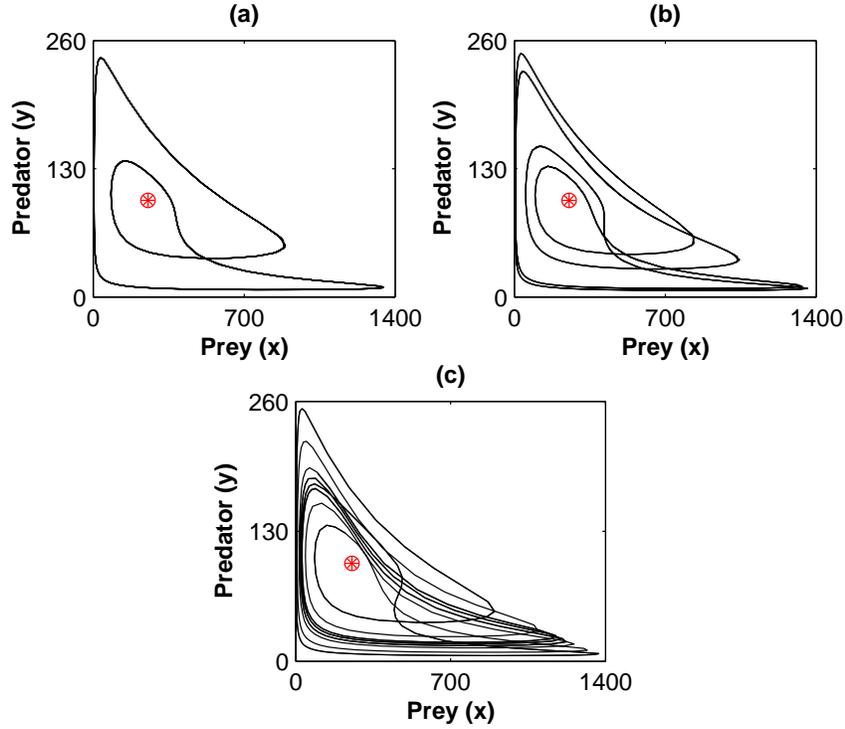}
\caption{nT periodic solution of the system $(1.3)$ with
different $\tau_2 ~(>\tau_{2_{0}})$ but fixed $\tau_1 ~(>\tau_{1_{0}})$.
(a) 2T period:
$\tau_1=0.5, \tau_2=0.56$,
(b) 4T period: $\tau_1=0.5, \tau_2=0.6$ and (c) 8T period: $\tau_1=0.5, \tau_2=0.6$. Other parameters
are as in Table 1.}
\label{fig2}
\end{figure}
 Our simulation results indicate that the system exhibits irregular periodic behaviors for
$\tau_{1}=0.7 ~(>\tau_{1_{0}})$ and $\tau_2=0.8 ~(>\tau_{1_{0}})$ (Fig. 4).
In order to characterize this irregular behavior, we perform the standard numerical diagnostics of the solutions, viz., Lyapunov exponents.
Lyapunov exponent or Lyapunov characteristic exponent of a dynamical system is a quantity that characterizes the rate of separation of infinitesimally close trajectories. Negative, zero and positive Lyapunov exponents ($\lambda$) indicate, respectively, the stable, unstable and chaotic behavior of the system. For the parameter values as in the Fig. 4, the value of the largest Lyapunov exponent ($\lambda$) is found
to be positive when $\tau_2$ varies, indicating chaotic
dynamics of the system (1.3) (Fig. 5).

%\noindent {\bf (ii) Sensitivity to initial conditions}
%
%Sensitivity to initial conditions means that arbitrarily closed points eventually separate from each other if the system is chaotic. Thus, an arbitrarily small perturbation of the current trajectory may lead to significantly different future behavior. To show the sensitivity of the system (1.3) to its initial value, we compute the error $\Delta x(t)=x_1(t)-x_2(t)$, where $x_1(t)= (30, 5.83)$ and  $x_2(t)= (30.01, 5.83)$ and plot them over time (Fig. 5(b)). This figure shows that initial error between two trajectories is negligible within the time interval $0<t<80$, but the error increases in the subsequent time. This phenomena ensures the chaotic nature of the system (1.3).\\
\begin{figure}
\centering
\includegraphics[width=5in,height=3in]{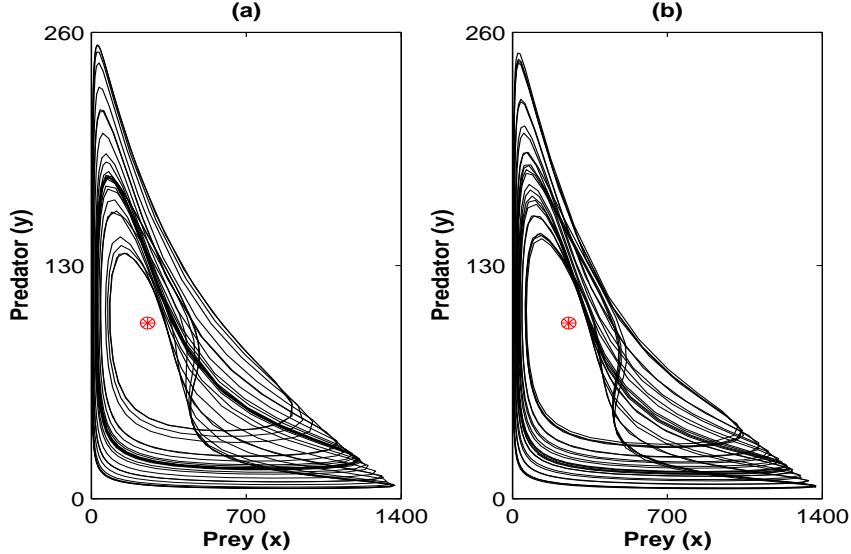}
\caption{The system $(1.3)$ exhibits chaotic dynamics for $\tau_2=0.626$ (Fig. (a)) and $\tau_2=0.66$ (Fig. (b)).
Other parameters are as in the Fig. 9.}
\label{fig2}
\end{figure}
%
%\noindent {\bf (iii) Power spectral density}
%
%Power spectral density (PSD) describes how the power of a signal or time series is distributed with frequency. The power spectra of the prey population is presented in Fig. 5(c). The irregular broad peaks of this figure are indicative of chaos and randomness.\\
%
%\noindent {\bf (iv) Recurrence plot}
%
%A recurrence plot (RP) is a plot showing the times at which trajectories visit roughly the same place in the phase space. In other words, it is a graph of $\overrightarrow{x}(i)\approx \overrightarrow{x}(j)$, showing $i$ on a horizontal axis and $j$ on a vertical axis, where $\overrightarrow{x}$ is a phase space trajectory. When the system is stable, then $\overrightarrow{x}(i) \approx \overrightarrow{x}(j)$ occurs for all succeeding $i$'s. If the system undergoes oscillatory state then $\overrightarrow{x}(i) \approx \overrightarrow{x}(j)$ happens in a particular rhythm whatever be the period of oscillations. In the chaotic scenario, $\overrightarrow{x}(i) \approx \overrightarrow{x}(j)$ occurs at arbitrary instants, $i.e.,$ the time intermission for nearly same system solutions does not follow any
%particular rule. The recurrence plot of the system for the parameter values as in the Fig. 5 is represented in Fig. 5d. To eliminate the impact of transients, we have discarded first 500 iterations. The random points on the Time-Time plane ensures that more or less same values of phase trajectories take place without rhythm, indicating the chaotic nature of the system.\\
\begin{figure}
\centering
\includegraphics[width=5in,height=3.5in]{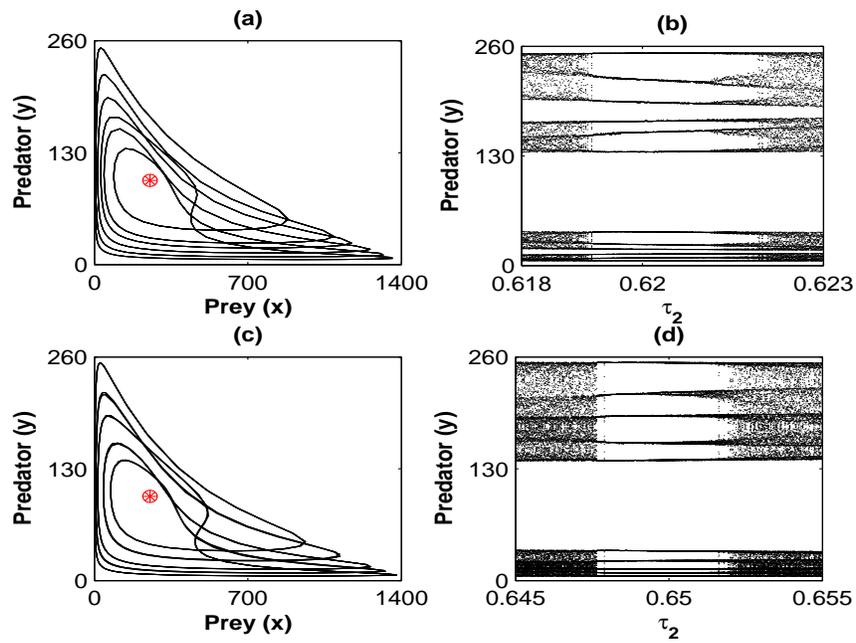}
\caption{The system $(1.3)$ exhibits $6T-$period solutions for $\tau_2=0.62$ (Fig. (a)) and $5T-$period solutions for $\tau_2=0.65$ (Fig. (c)). Magnified parts of the Fig. 6 over [0.618, 0.623] (Fig. (b)) and [0.645, 0.655] (Fig. (d)). Other parameters are as in the Fig. 4.}
\label{fig2}
\end{figure}
In Fig. 6, we have shown a typical bifurcation diagram of the predator population $y$ for fixed $\tau_1=0.5$ and variable $\tau_2$. The bifurcation diagram shows that the system (1.3) exhibits a rich dynamics such as period-doubling bifurcation, period-halving bifurcation, chaotic band, narrow windows, etc. as the parameter $\tau_2$ is smoothly varied. To observe the interplay between the strength of delay parameters, we have plotted the stability region of the system $(1.3)$ in the $\tau_2\tau_1$-plane for a fixed value of $m=0.45$ (Fig. 7). This figure indicates that the two delays follow inverse relationship for maintaining the stability of the system. However, if the delays are large then the system is always chaotic. Different dynamical features like stable equilibrium, different $nT$-periodic solutions and chaotic dynamics are clearly depicted by different indices in this figure.
The phase-space of the system (1.3) for different $\tau_2 ~(>\tau_{2_{0}})$ but fixed $\tau_1 ~(>\tau_{1_{0}})$ shows that the system exhibits $T$-period solution for $\tau_2=0.3$, $2T$-period solution for $\tau_2=0.4$ and $4T$-period solution for $\tau_2=0.53$ (Fig. 7). The system shows period-halving for $\tau_2=0.56$ and again period-doubling at $\tau_2=0.6$ and $\tau_2=0.626$ (Fig. 8). The system exhibits chaotic dynamics for $\tau_2=0.63$ and $\tau_2=0.66$ (Fig. 9). The system again shows $6T$-period solution for $\tau_2=0.62$ (Fig. 10(a)) and $5T$-period solution for $\tau_2=0.65$ (Fig. 10(c)). Periodic windows are intermittently scattered. The magnified periodic windows over the values of $\tau_2= [0.618, 0.623]$ and $\tau_2= [0.645, 0.655]$ are shown in Fig. 10(b) and Fig. 10(d), respectively.\\

%To observe the inter play between the strength of prey refuge and delay parameters, we have plotted the stability region of the system $(1.3)$ in the $m\tau_1$-plane when $\tau_2$ is fixed (Fig. 11(a)). In Fig. 11(b), the same is plotted in $m\tau_2$-plane when $\tau_1$ is fixed. These figures show that the value of the strength of prey refuge increases with the length of one delay when other delay is fixed. Fig. 11(c) gives the stability region of the system in the $\tau_2\tau_1$-plane for fixed value of $c$. This figure indicates that the two delays follow inverse relationship for maintaining the stability of the system. However, if the delays are large then the system is always unstable.

%\begin{figure}
%\centering
%\includegraphics[width=5in,height=3.5in]{./fig13.eps}
%\caption{Stability regions of the system $(1.3)$ are depicted in different parametric planes. Fig. (a) is the stability region in the $m\tau_1$-plane for fixed $\tau_2 ~(=0.18$). Fig. (b) is the stability region in the $m\tau_2$-plane for fixed $\tau_1 ~(=0.24)$.
%Fig. (c) is the stability region in the $\tau_1\tau_2$-plane for  fixed $m ~(=0.45)$.
%Remaining parameters are as in the Table 1.}
%\label{fig2}
%\end{figure}

\section{Summary}

Effect of prey refuge has not been considered explicitly in the prey-predator
models. However, both field and laboratory experiments confirm that prey refuge reduces predation rates by decreasing encounter rates between predator and prey. On the other hand, a prey-predator model becomes more realistic in the presence of different
delays which are unavoidable elements in physiological and ecological processes. In this
paper, a prey-predator model that incorporates different biological delays and the effect
of prey refuge is studied. A time delay $\tau_1$ is considered in the logistic prey growth
rate to represent density dependent feedback mechanism and the second time delay $\tau_2$
is considered in the predator response function to represent its gestation delay. The
objective is to study the dynamic behavior of a multi-delayed prey-predator system in the
presence of prey refuge.

It is observed that the non-delayed system is asymptotically stable under some parametric restrictions. There is a critical value ($\tau_{2_0}$) of the gestation delay parameter ($\tau_2$), below which the single-delayed system ($\tau_1=0$) is locally asymptotically stable and above
which the system is unstable. A Hopf-bifurcation occurs when the delay parameter attains
the critical value $\tau_{2_0}$. Keeping gestation delay ($\tau_2$) within its stability range, a critical
value $\tau_{1_0}$ of the negative feedback mechanism delay parameter ($\tau_1$) is obtained below
which the double-delayed system is locally stable and above which the system is unstable.
A Hopf-bifurcation occurs at $\tau_1=\tau_{1_0}$. Similarly, if we first keep $\tau_2=0$, we get previous type result from this system. The system exhibits irregular behavior when these delays are large and above their critical values. This irregularity has been identified as
chaotic through different tests. These simulations indicate that the system exhibits a rich
dynamics such as period-doubling bifurcation, period-halving bifurcation, chaotic band,
narrow and wide windows etc as the parameter $\tau_2$ is smoothly varied. It is also observed
that the strength of the prey refuge increases with delays to keep the system in
stable condition. The inter-play between two delays for fixed value of prey refuge has also been determined. It is noticed that these delays work in a complementary fashion. In other words, to keep the system in stable condition, the delay in the logistic
prey growth should be low when gestation delay is high or vice versa. Thus, a prey-predator system may exhibit simple stable behavior, regular cyclic behavior
or chaotic behavior depending on the length of delays.

\newpage

\newpage
{\bf Acknowledgment}\\
M. Lakshmanan acknowledges National Academy of Sciences, India for financial support in the form of a Senior Scientist Platinum Jubilee Fellowship.

%\newpage
%\vspace{.5cm}
%\begin{center} {\bf Table 1.} {\it Equilibrium points of the model
%system (2.3) and their existence conditions}.\end{center}
%
%$$ \begin{array}{ccccc} \hline
%Equilibrium~ Points & Coordinates & Conditions~ for~ Existence \\ \hline
%  E_0 & (0,0) & always~ exists \\
%  E_1 & (k,0) & always~ exists \\
%  E^* & (x^{*},y^{*}), where & exists~ if~ c < c_1,   \\
%   & x^{*}=\frac{d }{\alpha(1-m)(\theta-hd)}, &  where~ c_1=1- \frac{d}{\alpha k(\theta-hd)},\\
% & y^{*}=\frac{r(k-x^*)\{1+\alpha h(1-m)x^*\}}{\alpha k(1-m)}.&\theta>hd+\frac{d}{\alpha k}&\\ \hline
%\end{array}$$

\newpage

%\begin{center} {\bf Table 2.} {\it Equilibrium points of the model
%system (2.3) and their stability conditions}.\end{center}
%
%$$ \begin{array}{ccccc} \hline
% Equilibrium~ points & Nature~ of~ equilibrium~ point & Stability ~Conditions \\ \hline
%  E_0 & unstable~ saddle & unconditionally  \\
%  E_1 & asymptotically~ stable & c > c_1\\
%  & &\\
%   & asymptotically~ stable & c_2< c < c_1   \\
%  E^* & unique~ limit~ cycle & 0 < c < c_2   \\
%  & Hopf~ bifurcation &  c = c_2,   \\\hline
%Hhere ~c_1=1- \frac{d}{\alpha k(\theta-hd)}, & c_2=1-
%\frac{\theta+hd}{\alpha kh(\theta-hd)}, &  &\\
% ~~~~~\theta>max[\frac{hd(\alpha kh+1)}{\alpha
%kh-1},hd+\frac{d}{\alpha k}] & and ~~ \alpha>\frac{1}{kh}. &  &\\  \hline
%\end{array}$$

%\vspace{1in}
%
%\begin{center} {\bf Table 3.} {\it  Parameter values used in the model
%simulation}.\end{center}
%
%$$ \begin{array}{ccccc} \hline
% Parameters & Values \\ \hline
%  r &  2.65  \\
%  k & 898 \\
%  \alpha & 0.045 \\
%  c & .... \\
%  h & 0.0437 \\
%  \theta & 0.215 \\
%  d & 1.06 \\
% & & \\  \hline
%\end{array}$$

\end{document}